\documentclass[letterpaper]{article}%
\usepackage{amssymb}
\usepackage{amsmath}
\usepackage{amsfonts}
\usepackage{sw20bams}
\usepackage{comment}
\usepackage{hyperref}
\usepackage{url}
\usepackage{aurl}%
\usepackage{graphicx}
\providecommand{\U}[1]{\protect\rule{.1in}{.1in}}

\ifx\pdfoutput\relax\let\pdfoutput=\undefined\fi
\newcount\msipdfoutput
\ifx\pdfoutput\undefined\else
\ifcase\pdfoutput\else
\ifx\paperwidth\undefined\else
\ifdim\paperheight=0pt\relax\else\pdfpageheight\paperheight\fi
\ifdim\paperwidth=0pt\relax\else\pdfpagewidth\paperwidth\fi
\fi\fi\fi
\begin{document}

\title{Fast Computing Formulas for some Dirichlet L-Series}
\author{Jorge Zuniga\thanks{\hspace{4pt}Independent Researcher}}
\date{\hspace{85pt}February 2026.\thanks{\hspace{4pt}Keywords: $L$--functions, Wilf--Zeilberger method, Series acceleration, Algorithms, Binary Splitting}\vspace{8pt} \newline Dedicated to Jesús Guillera, 31\thinspace Dec\thinspace 1955\thinspace--\thinspace9\thinspace Feb\thinspace2026}
\maketitle

\begin{abstract}
This work introduces new hypergeometric-type rational series for fast computing of some primitive Dirichlet \textit{L}--Series formulas, $L_k(s):=L(s,\chi_k)$ for $s=2,3$ using the binary splitting algorithm. Series were found and proven using the WZ--Dougall's $_5H_5$ method to accelerate Hurwitz $\zeta$ identities. The new algorithms include $\zeta(3):=L_1(3)$ (Ap\'{e}ry's constant), $G:=L_{-4}(2)$ (Catalan's constant) as well as $L_{-7}(2)$, $L_{-8}(2)$, $L_{-15}(2)$, $L_{-20}(2)$, $L_{-24}(2)$ together with $L_5(3)$, $L_8(3)$ and $L_{12}(3)$. Formulas were tested and verified up to 100 million decimal places for each $L$-value.

\end{abstract}

\section{Introduction\smallskip}
\vspace*{4pt}
Jes\'{u}s Guillera \cite{GUI} has recently published a method to accelerate the Hurwitz zeta function $\zeta(s,x)$ for $s=2,3$ and $x\in\mathbb{Q}^+$ with Wilf Zeilberger \cite{WILFZEIL}%
\cite{KOEPF} pairs $(F,G):=(F(n.k),G(n,k))$ that are transformed via $\text{WZ}_{u,v}:F(n,k)\rightarrow F(u\,n,k+v\,n)$ for some $u\in\mathbb{N}$ and $v\in\mathbb{Z}$. For selected $(u,v)$ the logarithmically convergent slow  series that defines $\zeta(s,x)$, Eq.(\ref{8z}), can be transformed into a new proven accelerated series having linear convergence. \vspace{+4pt}

In this work another approach is taken to deal with Hurwitz $\zeta$'s aiming to discover new strongly convergent identities for self--dual primitive Dirichlet \textit{L}--Series, $L_k(s):=L(s,\chi_k)$, $s=2,3$. The WZ method is instead applied on a partial Dougall's $_5H_5$ hypergeometric parametrized sum to prove some of the most efficient $L$--Series formulas known so far. This $_5H_5$ approach was first applied by Chu
and Zhang \cite{CHUZANG} using Abel's summation by parts to find and prove many hypergeometric identities. Guillera \cite{GUILLERA} merged it later with WZ techniques by building an automated tool to find and prove new identities that is a powerful resource to work with WZ pairs overrunning the $\text{WZ}_{u,v}$ acceleration methods initially introduced by Amdeberhan and
Zeilberger \cite{ZEIL}. The WZ--Dougall's $_5H_5$ procedure has been recently extended by Campbell \cite{CAMP1}\cite{CAMP2}. \vspace{+4pt} 

The interested reader is encouraged to download and install Guillera's Maple code\daurl{:(here)}{https://anamat.unizar.es/jguillera/Maple/WZ-DOUGALL-CONSTANTS-OCTOBER-2023.mw} \aurl{:(here)}{}, run the \texttt{SER} commands given further on to get the series presented in this work, check the proofs certificates (there is no room to put them here) and analyze the $(F,G)$ details. y--cruncher's \cite{YEE0} configuration files are also provided to test the new formulas and compute for first time these $L$--function values to several millions of decimal places. See the \textit{Appendix} for the association of these $L$--values to the LMF Data Base \cite{LMFDB}\vspace*{-2pt}

\section{Hypergeometric Sums\smallskip}
\vspace*{4pt}
Let $p(n)$, $q(n)$ and $r(n)$ be integer coefficient polynomials such that
none of them vanishes for$\ n\in 
\mathbb{N}
$, and polynomials $r(n)$ and $q(n)$ are of the same degree ${\small \textit{D}}$ having
rational roots $1-r_{\ell}$ and $1-q_{m}$ respectively with $0<r_{\ell} 
,q_{m}\leq1$ and $r_{\ell}\neq q_{m}$ for $\ell,m=1,2,...,{\small \textit{D}},$ so that
\{$r_{1}$,$r_{2}$...,$r_{\tiny \textit{D}}$\}\ and\textsf{\ }\{$q_{1}$,$q_{2}$...,$q_{\tiny \textit{D}}$\}
are disjoint rational multisets.\vspace{+6pt} 

With this setup the real series $\mathcal{S}
=\lim_{N\rightarrow\infty}\mathcal{S}_{N},$ is defined by convention through the following
\textit{rational} approximants sequence of hypergeometric-type convergent
sums $\mathcal{S}_{N}$\vspace*{4pt}
\begin{equation}
	\mathcal{S}_{N}=\frac{1}{\beta_0}\cdot
	{\displaystyle\sum\limits_{n=1}^{N}}
	\mathcal{R(}n\mathcal{)\cdot H(}n)\text{ }\mathcal{=}\text{ }\frac{1}{\beta_0
	}\cdot%
	{\displaystyle\sum\limits_{n=1}^{N}}
	\mathcal{R(}n\mathcal{)\cdot\rho}^{n}\,\mathcal{M}_n,\text{
		\ \ {\small \textit{N}}}\in%
	\mathbb{N}%
	\label{1}%
\end{equation}
where $\beta_0\in%
\mathbb{Z}
_{\neq0}$ is a generic normalization factor, $\mathcal{R(}n\mathcal{)}%
=p(n)/r(n)$ is the rational part\ and $\mathcal{H(}n\mathcal{)}=\mathcal{\rho
}^{n}\,\mathcal{M}_n$ is the hypergeometric factor where $|\rho
|\leq1,$ $\rho$ $\in%
\mathbb{Q}
$ is the convergence rate and $\mathcal{M}_n$ is defined by the
following ratio of products of Pochhammer's symbols$\vspace*{4pt}$%
\begin{equation}
	\mathcal{M}_n=\left[
	\begin{array}
		[c]{cccc}%
		r_{1} & r_{2} & ... & r_{\small \textit{D}}\\
		q_{1} & q_{2} & ... & q_{\small \textit{D}}%
	\end{array}
	\vspace*{1pt}\right]  _{n}=\dfrac{(r_{1})_{n}\,(r_{2})_{n}\,...\,(r_{\small \textit{D}})_{n}%
	}{(q_{1})_{n}\,(q_{2})_{n}\,...\,(q_{\small \textit{D}})_{n}}\vspace*{6pt}\label{2}%
\end{equation}
where ${\small \textit{D}}$ is the hypergeometric deepness and $(\upsilon)_{n}=\upsilon
(\upsilon+1)\ldots(\upsilon+n-1)$ is the rising factorial. This means
that\ there exist leading coefficients $c_{r},c_{q}\in%
\mathbb{Z}
_{\neq0}$ whose ratio holds $c_{r}/c_{q}=\rho=\lim_{n\rightarrow\infty
}r(n)/q(n)$ and\vspace*{-2pt}%

\begin{equation}%
	\begin{array}
		[c]{ccc}%
		r(n) & = & c_{r}(n-1+r_{1})(n-1+r_{2})\ldots(n-1+r_{\small \textit{D}})%
		\genfrac{}{}{0pt}{0}{\genfrac{}{}{0pt}{1}{{}}{{}}}{{}}%
		\\
		q(n) & = & c_{q}(n-1+q_{1})(n-1+q_{2})\ldots(n-1+q_{\small \textit{D}})%
		\genfrac{}{}{0pt}{0}{\genfrac{}{}{0pt}{1}{{}}{{}}}{{}}%
	\end{array}
	\vspace*{2pt}\label{2a}%
\end{equation}
Note that rational roots $1-r_{\ell}$ and $1-q_{m}$ enter as elements
$r_{\ell}$ and $q_{m}$ in $\mathcal{M}_n$, therefore$\vspace*{2pt}$%
\begin{equation}
	\mathcal{H(}n\mathcal{)}=\rho^{n}\cdot\dfrac{(r_{1})_{n}\,(r_{2}%
		)_{n}\,...\,(r_{\small \textit{D}})_{n}}{(q_{1})_{n}\,(q_{2})_{n}\,...\,(q_{\small \textit{D}})_{n}}=%
	{\displaystyle\prod\limits_{k=1}^{n}}
	\frac{r(k)}{q(k)}\vspace*{2pt}\label{2b}%
\end{equation}
so the hypergeometric-type series can be equivalently written either
as$\vspace*{2pt}$%
\begin{equation}
	\mathcal{S}=\frac{1}{\beta_0}\cdot%
	{\displaystyle\sum\limits_{n=1}^{\infty}}
	\frac{p(n)}{r(n)}\cdot\mathcal{\rho}^{n}\mathcal{\cdot}\left[
	\begin{array}
		[c]{cccc}%
		r_{1} & r_{2} & ... & r_{\small \textit{D}}\\
		q_{1} & q_{2} & ... & q_{\small \textit{D}}%
	\end{array}
	\vspace*{1pt}\right]  _{n}\vspace*{3pt}\label{7a}%
\end{equation}
or\vspace{-4pt}%
\begin{equation}
	\mathcal{S}=\frac{1}{\beta_0}\cdot%
	{\displaystyle\sum\limits_{n=1}^{\infty}}
	\frac{p(n)}{r(n)}\cdot%
	{\displaystyle\prod\limits_{k=1}^{n}}
	\frac{r(k)}{q(k)}\vspace*{3pt}\label{7b}%
\end{equation}
This last identity is the classical form to apply the Binary Splitting algorithm
\cite{PAPA} and variations \cite{CHENG} where the computing time to evaluate
series $\mathcal{S}$ with $B$-bit is $\mathcal{O}(\mathit{M}(B)\cdot\log B)$,
being \textit{M}$(B)$ the time needed to multiply two $B$-bit numbers.
\subsection{Binary Splitting Cost\vspace{0pt}}
 A practical measurement associated to \textit{M}$(B)$, that allows hypergeometric series
 of different types with the same run-time complexity to be directly compared and
 ranked by performance is\ given by the \textit{Binary Splitting Cost
}\cite{AYEE}\vspace{-2pt}
\begin{equation}
	C_{s}=-\frac{4\cdot {\small \textit{D}}}{\log\rho}\label{3}%
\end{equation}
$C_{s}$ provides a method to know a priori which algorithm (hypergeometric identity, formula)
is better or faster than another, being a fundamental quantity to search, find
and select the most efficient linearly convergent hypergeometric-type rational series from a pool of possibilities and options.
\pagebreak
\section{Dougall's $_5H_5$ and the WZ Method\smallskip}

Dougall's \({}_{5}H_{5}\) identity (often referred to as Dougall’s bilateral sum \cite{DOUG}) evaluates a specific very well-poised hypergeometric series in terms of a ratio of products of gamma functions (see \cite{CHUZANG}\cite{GUILLERA} for notation
and references therein). The truncated unilateral sum, depending on five free parameters, is also a very well poised series and serves as a prolific source to generate new proven identities producing hypergeometric accelerations via WZ as follows,  
\begin{equation}
	\Omega(a,b,c,d,e)=\sum_{k=0}^{\infty}\Omega_k(a,b,c,d,e)\label{5h1}%
\end{equation}\vspace{0pt}
\begin{equation}
	\Omega_k(a,b,c,d,e)=\frac{(b)_{k}(c)_{k}(d)_{k}(e)_{k}}{(1+a-b)_{k}(1+a-c)_{k}(1+a-d)_{k}%
		(1+a-e)_{k}}(a+2k)\label{5h2}%
\end{equation}
with
\begin{equation}
\mathfrak{R}(1+2a-b-c-d-e)>0\label{5}%
\end{equation}
for convergence. The $\Omega(a,b,c,d,e)$ hypergeometric series is symmetric in $\left(b,c,d,e\right)$ 
and can be expressed in terms of another shifted hypergeometric-type series of the same
kind
\begin{equation}
	\Omega(a+\alpha\,n,b+\beta\,n,c+\gamma\,n,d+\delta\,n,e+\varepsilon\,n)\label{5o}%
\end{equation} 
where the iteration pattern $(\alpha,\beta,\gamma,\delta,\varepsilon)$
just contains non-negative integers. For $e=a$ the particular case of Dougall’s \({}_{5}F_{4}\) four parameters monotone series identity (also known as Dougall’s theorem \cite{DOUG}\cite{5F4}) is obtained, while for $e\rightarrow-\infty$ the special
sum $\Omega(a,b,c,d,-\infty)$ is alternating for real finite
positive parameters. In this case%
\begin{equation}
\Omega_k(a,b,c,d,-\infty)=
\frac{(-1)^{k}(b)_{k}(c)_{k}(d)_{k}}{(1+a-b)_{k}(1+a-c)_{k}(1+a-d)_{k}%
}(a+2k)\label{5a}%
\end{equation}
By applying the Wilf Zeilberger method to the resulting shifted summands in Eq.(\ref{5h2}) and
Eq.(\ref{5a}) --see \cite{GUILLERA} for details-- with pairs $(F,G):=(\hspace{1pt}F(n,k),G(n,k)\hspace{1pt})$ and taking\footnote{Notice that Eqs.(\ref{5g}--\ref{7}) are closed under $\text{WZ}_{u,v}:F(n,k)\rightarrow F(u\hspace{1pt}n,k+v\hspace{1pt}n)$ transformation. Pattern $(\alpha,\beta,\gamma,\delta,\varepsilon)\in\mathbb{N}^5$ is shifted to $(\alpha,\beta,\gamma,\delta,\varepsilon)\cdot u +(2,1,1,1,1)\cdot v$ for $u\in\mathbb{N}$ and $v\in\mathbb{Z}_{\vspace{-6pt}{>-u\cdot \min(\alpha/2,\beta,\gamma,\delta,\varepsilon)}}$.}
\begin{equation}
	F(n,k)=\Omega_k(a+\alpha\,n,b+\beta\,n,c+\gamma\,n,d+\delta\,n,e+\varepsilon\,n)\label{5g}%
\end{equation}
only recurrences of first order are obtained giving pairs $(F,G)$ together with proven hypergeometric series identities $\Omega(a,b,c,d,e)=\sum_{n=0}^\infty G(n,0)$ being $G(n,0)$ of the type of Eqs.(\ref{7a}--\ref{7b}). The whole process is fully automated and the binary splitting costs $C_{s}$ are calculated to classify and select the discovered series. Maple\texttrademark\ program commands are, 
\begin{equation} \label{6}
	\begin{array}{l}
		\texttt{SER(}a + \alpha n, b + \beta n, c + \gamma n, d + \delta n, e + \varepsilon n\texttt{)} \\
		\\
		\texttt{SER(}a + \alpha n, b + \beta n, c + \gamma n, d + \delta n, \texttt{oo)}    \text{ if } e \to -\infty
	\end{array}
\end{equation}
respectively. Note that in order to keep the convergence of the
resulting series, the iteration pattern
must hold for base parameters satisfying Eq.(\ref{5}), either\vspace*{-2pt}%
\begin{equation}%
\genfrac{}{}{0pt}{0}{\genfrac{}{}{0pt}{1}{{}}{{}}}{{}}%
\alpha>\tfrac{1}{2}(\beta+\gamma+\delta+\varepsilon)\text{ \ \ or\ \ \ }%
\alpha>\tfrac{2}{3}(\beta+\gamma+\delta)\vspace*{-2pt}\label{7}%
\end{equation}
for Eq.(\ref{5h2}) or Eq.(\ref{5a}) respectively. For given base parameters $(a,b,c,d,e)$ each series that matches Eqs.(\ref{5h1}--\ref{5h2}) or Eq.(\ref{5a}) can be strongly accelerated in this way by finding the optimal searching pattern $(\alpha,\beta,\gamma,\delta,\varepsilon)$ that gives the minimum cost. In what follows, this general method is applied on some Hurwitz zetas to get fast Dirichlet \textit{L}--Series formulas\vspace{4pt}.

\section{Computing the Hurwitz Zeta Function  \bigskip} The Hurwitz zeta function is traditionally 
defined as
\begin{equation}
	\zeta(s,x)=\sum_{k=0}^\infty \frac1{(k+x)^s},\ \ \ \ \Re(s)>1\ \ \ x\notin\mathbb{Z}_{<1}\label{8z}%
\end{equation}
Series converges absolutely for $\Re(s) = \sigma > 1$ and, if $\delta>0$, uniformly in every
half-plane $\sigma\ge1+\delta$, so that $\zeta(s,x)$ is an analytic function of $s$ in the half-plane
$\Re(s) > 1$. 
\subsection{Reduction Formulas\vspace{0pt}}
The Hurwitz zeta function holds a reflection identity in the variable $x$ for $s\in\mathbb{N}$ that is a direct consequence of the reflection formula for the polygamma function,\vspace{-4pt}  
\begin{equation}
\zeta(s,1-x)+(-1)^s\zeta(s,x)=\pi\cdot\frac{(-1)^{s-1}}{(s-1)!}\cdot\frac{d^{s-1}}{dx^{s-1}}\cot(\pi x)\vspace{2pt}\label{8r}%
\end{equation}
for $s=2$ a mirror relationship for $\mathcal{S}_0(x,y):=\zeta(2,x)-\zeta(2,y)$ is obtained from\vspace{2pt}
\begin{equation}
	\zeta(2,1-x)+\zeta(2,x)=\pi^2\cdot\csc^2(\pi x)\vspace{2pt}\label{8s}%
\end{equation}
\begin{equation} 
	\mathcal{S}_0(x,y)-\mathcal{S}_0(1-y,1-x)=\pi^2\cdot\left(\csc^2(\pi x)-\csc^2(\pi y)\right)\vspace{4pt}\label{8w}%
\end{equation}
for $s=3$ 
\begin{equation}\label{8t}
	\zeta(3,1-x)-\zeta(3,x)=\pi^3\cdot\cot(\pi x)\cdot\csc^2(\pi x)\vspace{2pt}
\end{equation}
where the rhs can be converted to radicals for some special rational values of $x,y\in(0,1)$. 
Hurwitz zeta also holds a multiplication identity that is derived from Gauss multiplication
 formula for the gamma function and its variation for polygammas. For $k\in\mathbb{N}$ and $s\ne1$  
\begin{equation}\label{8u}
	\zeta(s,k x)=k^{-s}\cdot\sum_{n=0}^{k-1}\zeta\left(s,x+\small{\frac{n}{k}}\right)\vspace*{4pt}
\end{equation}
By setting $x=k^{-1}$ 
\begin{equation}\label{8v}
	(k^s-1)\zeta(s)=\sum_{n=1}^{k-1}\zeta\left(s,\small{\frac{n}{k}}\right)\vspace*{4pt}
\end{equation}
where $\zeta(s):=\zeta(s,1)$ is Riemann's zeta function. Eqs.(\ref{8u}--\ref{8v}) are some kind of reduction 
formulas as well since they allow to express some sums of Hurwitz zetas with $0<x\le 1$ and $x\in\mathbb{Q}$ related to lower and faster 
computing constants like $\pi$, Apery's constant $\zeta(3)$, Catalan's constant $G$ and $L_{-3}(2)$ 
as it will be seen further on. 
\subsection{Acceleration Seeds\vspace{0pt}} For $s=2$ and $s=3$, Hurwitz $\zeta(s,x)$ Eq.(\ref{8z}) matches Dougall Eqs.(\ref{5h1}--\ref{5h2}) producing these acceleration seeds for $0< x,y\le 1$ and $x,y\in\mathbb{Q}$
\begin{equation}
	\zeta(2,x)-\zeta(2,y)=\frac{y-x}{x^2\,y^2}\cdot\Omega(x+y,x,x,y,y)\label{8xy}%
\end{equation}
In this case the difference of two Hurwitz zetas yields one accelerated series. This is also a reduction formula. For $s=3$ there is not reduction, a single accelerated series is obtained
\begin{equation}\vspace{3pt}
	\zeta(3,x)=\frac{1}{2\,x^4}\cdot\Omega(2x,x,x,x,x)\label{8x}%
\end{equation}Notice that rational factors on the rhs can be ignored or omitted to search the accelerated series since they are absorbed by $p(n)$ and $\beta_0$ in Eqs.(\ref{7a}--\ref{7b}). They can be restored later.\pagebreak
\subsection{$\boldsymbol{\mathcal{S}_0(x,y)=\zeta(2,x)-\zeta(2,y)}$\vspace{0pt}}From here on out $\mathcal{S}_0(x,y)$ will be referred to the specific series generated by pattern $(\alpha,\beta,\gamma,\delta,\varepsilon)=(3,1,1,1,1)$. In this case this command is applied 
\begin{equation}\label{9x}
\texttt{SER(}x + y + 3 n, x + n, x + n, y + n, y + n\texttt{)} 
\end{equation}
the following proven accelerated series is obtained
\begin{equation}\label{9y}
	\mathcal{S}_0(x,y)=\frac{y-x}{2\,x^2\,y^2}\cdot\sum_{n=0}^{\infty}\,\frac{\mathcal{P}(n,x,y)\cdot\mathcal{M}_n(x,y)}{(2n+1)(2n+x+1)^2(2n+y+1)^2}\left(-\frac{1}{2^{10}}\right)^n
\end{equation}\vspace{-8pt}
\begin{align}\label{z02}
\mathcal{M}_n(x,y)=\left[\begin{array}
		[l]{ccccccccc}
		1+x-y&1+y-x&x&x&y&y&1&1&1\\
		\frac{1+x}2&\frac{1+x}2&\frac{2+x}2&\frac{2+x}2&\frac{1+y}2&\frac{1+y}2&\frac{2+y}2&\frac{2+y}2&\frac12
	\end{array}\right]_n\vspace{0pt}
\end{align}\vspace{-8pt}
\begin{multline*}
	\mathcal{P}(n,x,y) = 205\,n^{6}+\left(287 x +287 y +496\right)\,n^{5}\\+\left(135 x^{2}+394 x y +135 y^{2}+579 x +579 y +485\right)\,n^{4}\\+\left(21 x^{3}+183 x^{2} y +183 x y^{2}+21 y^{3}+222 x^{2}+620 x y +222 y^{2}+458 x +458 y +240\right)\,n^{3}\\+(28 x^{3} y +85 x^{2} y^{2}+28 x y^{3}+27 x^{3}+215 x^{2} y +215 x y^{2}+27 y^{3}+137 x^{2}+356 x y +137 y^{2}\\+176 x +176 y +60)\,n^{2}\\+(13 x^{3} y^{2}+13 x^{2} y^{3}+22 x^{3} y +66 x^{2} y^{2}+22 x y^{3}+12 x^{3}+82 x^{2} y +82 x y^{2}+12 y^{3}\\+38 x^{2}+88 x y +38 y^{2}+32 x +32 y +6)\,n \\+2 x^{3} y^{3}+5 x^{3} y^{2}+5 x^{2} y^{3}+4 x^{3} y +13 x^{2} y^{2}+4 xy^{3}+2 x^{3}+10 x^{2} y +10 x y^{2}+2 y^{3}\\+4 x^{2}+8 x y +4 y^{2}+2 x +2 y
\end{multline*}

This series has hypergeometric deepness ${\small \textit{D}}=9$ (the number of columns in $\mathcal{M}_n$). It has a binary splitting cost Eq.(\ref{3}) $C_s=36/\log(1024)=5.1937$ for a wide range of values of $(x,y)$. However, for some particular $(x,y)\in\mathbb{Q}^2$ values, Pochhammers in numerator and denominator either cancel themselves or are passed to the summand's rational part as factors of the numerator's polynomial $\mathcal{P}$. In such cases ${\small \textit{D}}$ reduces. For example ${\small \textit{D}}=8$ for Catalan's constant $G=L_{-4}(2)$ (See \cite{HESS} Eq.(13), Pilehrood long series) and also for some new $L_{-8}(2)$ and $L_{-20}(2)$ formulas introduced further on. In these cases series cost goes down to $C_s=32/\log(1024)=4.6166$. 
\subsection{$\boldsymbol{\mathcal{S}_1(x)=\zeta(3,x)}$\vspace{0pt}}Selecting $(\alpha,\beta,\gamma,\delta,\varepsilon)=(4,2,1,1,1)$ this command is applied 
\begin{equation}\label{10x}
	\texttt{SER(}2 x + 4 n, x + 2 n, x + n, x + n, x + n\texttt{)} 
\end{equation}
and $\mathcal{S}_1(x)$ is defined hereinafter by this specific proven accelerated series for $\zeta(3,x)$
\begin{align}\label{10y}
\mathcal{S}_1(x)=\frac{1}{6\,x^3}\cdot\sum_{n=0}^{\infty}\,\frac{\mathcal{P}(n,x)\cdot\mathcal{M}_n(x)}{(3n+1)(3n+2)(3n+x+1)^3(3n+x+2)^3}\left(-\frac{2^6}{3^{12}}\right)^n
\end{align}\vspace{0pt}	
\begin{equation}\label{10z}
	\mathcal{M}_n(x)=\left[\begin{array}
		[l]{ccccccccccc}%
		1 & 1 & 1 & 1 & 1 & \frac12 & \frac12 & \frac12 & x & x & x
		\genfrac{}{}{0pt}{1}{\genfrac{}{}{0pt}{0}{{}}{{}}}{{}}%
		\\
		\frac13 & \frac23 & \frac{1+x}3 & \frac{1+x}3 & \frac{1+x}3 &
		\frac{2+x}3 & \frac{2+x}3 & \frac{2+x}3 & \frac{3+x}{3} & \frac{3+x}3 & \frac{3+x}{3}%
		\genfrac{}{}{0pt}{1}{\genfrac{}{}{0pt}{0}{{}}{{}}}{{}}%
	\end{array}\right]_n
\end{equation}\vspace{-4pt}	
\begin{multline*}
\mathcal{P}(n,x) =40885\ n^{8}+\left(84915 x + 162086\right) n^{7}+\left(75717 x^{2}+293364 x +276772\right) n^{6}\\+\left(37699 x^{3}+223728 x^{2}+427080 x +265784\right) n^{5}\\+\left(11319 x^{4}+93158 x^{3}+270156 x^{2}+339489 x +156997\right) n^{4}\\+\left(2049 x^{5}+22662 x^{4}+90088 x^{3}+170406 x^{2}+159186 x +58436\right) n^{3}\\+\left(207 x^{6}+3168 x^{5}+16596 x^{4}+42512 x^{3}+59160 x^{2}+44088 x +13392\right) n^{2}\\+\left(9 x^{7}+228 x^{6}+1584 x^{5}+5250 x^{4}+9766 x^{3}+10716 x^{2}+6696 x +1728\right) n\\ 6 x^{7}+60 x^{6}+255 x^{5}+603 x^{4}+872 x^{3}+792 x^{2}+432 x + 96	
\end{multline*}

In this case ${\small \textit{D}}=11$ and cost is $C_s=44/\log(531441/64)=4.8756$ except for some special values of $x$. For example ${\small \textit{D}}=8$ for $x=1/2$ or $x=1$ which gives sums for $\zeta(3)$ that are somewhat faster. 

\section{Computing Dirichlet L-Series \bigskip} From the class of Dirichlet \textit{L}--functions defined by analytic continuation of the following series converging for $\Re(s)>1$ 
\begin{equation}\label{L1}
	L(s,\chi)=\sum_{n=1}^\infty\frac{\chi(n)}{n^s}
\end{equation}
where $\chi(n)$ is a Dirichlet character mod $k$, only self--dual (real) primitive \textit{L}--series with period $|k|$ for some selected moduli $k$, where $\chi(n):=\chi_k(n)=\left(\frac{k}{n}\right)$ is the Kroenecker symbol, are considered in this work. These \textit{L}--series are expressed in terms of Hurwitz zetas (Apostol \cite{APO}, Ch 12) as \vspace{-6pt}
\begin{equation}\label{L2}
	L_k(s):=L(s,\chi_k)=\frac1{|k|^s}\sum_{n=1}^{|k|-1}\left(\frac{k}{n}\right)\,\zeta\left(s,\frac{n}{|k|}\right),\ \ \ |k| > 1
\end{equation} and $L_1(s)=\zeta(s)$ for $k=1$. The \textit{Appendix} contains reduced identities (including mirror identities) for selected $L$--Series expressed as Eq.(\ref{L2}) where Hurwitz $\zeta$ relationships, Eq.(\ref{8s}) to Eq.(\ref{8v}), were applied. These identities feed the WZ--Dougall's $_5H_5$ process.
\subsection{Primary $\boldsymbol{L}$--Series}Before presenting the reduced identities for $L_k(s),\ s=2,3$, it is necessary to set up (some not well) known fast formulas for primary \textit{L}--Series:$\ \pi,\ \zeta(3)$, $L_{-3}(2)$ and Catalan's constant $G$ that take part in the algorithms.
\subsubsection{$\boldsymbol{L_{-4}(1)},\ $constant $\boldsymbol{\pi}$}$L_{-4}(1)$ is the \textit{L}--Series equivalent to Madhava--Leibniz $\frac{\pi}4$ series. $\pi$ is computed by hypergeometric Chudnovsky's modular algorithm \cite{CHUD} that has a very low binary splitting cost $C_s=12/\log(151931373056000)=0.36748...$

\subsubsection{$\boldsymbol{L_{1}(3)},$ Apery's constant $\boldsymbol{\zeta(3)}$}$\zeta(3)$ is computed by means of identities \textbf{\textit{i})} and \textbf{\textit{ii})} (next page) found by this author via WZ--Dougall $_5H_5$ in 2023. They are the fastest known series so far to calculate this constant. Both of them use the acceleration seed Eq.(\ref{8x}) with $x=1$. In fact some highly efficient classical formulas are also proven in this way by taking\footnote{For example $(\alpha,\beta,\gamma,\delta,\varepsilon)=(3,1,1,1,1)$ proves Amdeberhan--Zeilberger formula \cite{ZEIL} while $(5,2,2,1,1)$ gives Wedeniwski's $\zeta(3)$ series \cite{WED}.}\vspace{0pt}
\begin{equation*}
	\texttt{SER(}2 + \alpha\,n, 1 + \beta\,n, 1 + \gamma\,n, 1 + \delta\,n, 1 + \varepsilon\thinspace n\texttt{)} \vspace{4pt}
\end{equation*}   
 
 \textbf{\textit{i})} The currently fastest $\zeta(3)$ series is proven with pattern $(\alpha,\beta,\gamma,\delta,\varepsilon)=(11,4,3,2,1)$, the command\vspace{2pt} 
\begin{equation}\label{z1} 
	\texttt{SER(}2 + 11 n, 1 + 4 n, 1 + 3 n, 1 + 2 n, 1 + n\texttt{)} 
\end{equation}brings $\rho=-1/717445350000$, {\small \textit{D}}$\ =14$ and cost  $C_s=56/\log(717445350000)=2.0513..$
with a series giving about 11.86 decimal digits per term,
\begin{equation}
	\zeta(3)=
	{\displaystyle\sum\limits_{n=1}^{\infty}}\ \rho^n\cdot
	\frac{\,\mathcal{P}(n)}{\,\mathcal{R}(n)}\cdot\mathcal{M}_n\vspace{-4pt}\label{z2}%
\end{equation}
\begin{equation}\label{z21}
	\mathcal{M}_n=\left[\begin{array}
		[l]{cccccccccccccc}%
		1 & 1 & 1 & 1 & 1 & \frac12 & \frac12 & \frac12 & \frac13 & \frac23 & \frac14 & \frac34 & \frac16 & \frac56 
		\genfrac{}{}{0pt}{1}{\genfrac{}{}{0pt}{0}{{}}{{}}}{{}}%
		\\
		\frac19 & \frac29 & \frac49 & \frac59 & \frac79 & \frac89 & \frac{1}{10} &
		\frac{3}{10} & \frac{7}{10} & \frac{9}{10} & \frac1{12} & \frac5{12} & \frac7{12} & \frac{11}{12}
		\genfrac{}{}{0pt}{1}{\genfrac{}{}{0pt}{0}{{}}{{}}}{{}}%
	\end{array}\right]_n
\end{equation}and
\begin{equation}%
	\begin{array}
		[c]{cl}%
		\text{{\small \text{$\mathcal{P}$}(\textit{n}) \textit{\ }= }} & \text{\thinspace
			{\small 1565994397644288\thinspace}}{\small n}^{11}%
		\,\text{{\small --\thinspace6719460725627136\thinspace}}{\small n}%
		^{10}\,\text{{\small +\thinspace12632254526031264\thinspace}}{\small n}^{9}%
		\genfrac{}{}{0pt}{0}{\genfrac{}{}{0pt}{1}{{}}{{}}}{{}}%
		\\
		& \text{{\small --\thinspace13684352515879536\thinspace}}{\small n}%
		^{8}\,\text{{\small +\thinspace9451223531851808\thinspace}}{\small n}%
		^{7}\,\text{{\small --\thinspace4348596587040104\thinspace}}{\small n}^{6}\,%
		\genfrac{}{}{0pt}{0}{\genfrac{}{}{0pt}{1}{{}}{{}}}{{}}%
		\\
		& \text{{\small +\thinspace1352700034136826\thinspace}}{\small n}%
		^{5}\,\text{{\small --\thinspace282805786014979\thinspace}}{\small n}%
		^{4}\,\text{{\small +\thinspace38721705264979\thinspace}}{\small n}^{3}\,%
		\genfrac{}{}{0pt}{0}{\genfrac{}{}{0pt}{1}{{}}{{}}}{{}}%
		\\
		& \text{{\small --\thinspace3292502315430\thinspace}}{\small n}^{2}%
		\,\text{{\small +\thinspace156286859400\thinspace}}{\small n}%
		\,\,\text{{\small --\thinspace3143448000}}%
		\genfrac{}{}{0pt}{0}{\genfrac{}{}{0pt}{1}{{}}{{}}}{{}}%
		\\
		\text{{\small \text{$\mathcal{R}$}(\textit{n}) \textit{\ }= }} &
		\text{{\small --\thinspace48\thinspace\textit{n}}}^{5}\text{{\small (2\textit{n}--1)}}^{3}%
		\text{{\small (3\textit{n}--1)(3\textit{n}--2)(4\textit{n}--1)(4\textit{n}--3)(6\textit{n}--1)(6\textit{n}--5)}}%
		\genfrac{}{}{0pt}{0}{\genfrac{}{}{0pt}{1}{{}}{{}}}{{}}%
	\end{array}\label{z3}%
	\end{equation}\vspace{4pt}
	
	\textbf{\textit{ii})} A second $\zeta(3)$ fast series rises from $(\alpha,\beta,\gamma,\delta,\varepsilon)=(10,4,3,2,1)$ with command\vspace{0pt}
	\begin{equation}\label{z4} 
		\texttt{SER(}2 + 10 n, 1 + 4 n, 1 + 3 n, 1 + 2 n, 1 + n\texttt{)}\vspace{0pt} 
	\end{equation}In this case Eq.(\ref{z2}) holds with $\rho=2^{-16}\,3^{-12}$, {\small \textit{D}}$\ =14$ and cost $C_s=56/\log(2^{16}\,3^{12})=2.30702..$. 
\begin{equation}
	\mathcal{M}_n=\left[\begin{array}
		[l]{cccccccccccccc}%
		1 & 1 & 1 & 1 & 1 & \frac12 & \frac13 & \frac23 & \frac14 & \frac34 & \frac15 & \frac25 & \frac35 & \frac45 
		\genfrac{}{}{0pt}{1}{\genfrac{}{}{0pt}{0}{{}}{{}}}{{}}%
		\\
		 \frac18 & \frac38 & \frac58 & \frac78 & \frac19 & \frac29 & \frac49 & \frac59 & \frac79 & \frac89 & \frac{1}{10} &
		\frac{3}{10} & \frac{7}{10} & \frac{9}{10} 
		\genfrac{}{}{0pt}{1}{\genfrac{}{}{0pt}{0}{{}}{{}}}{{}}%
	\end{array}\right]_n\,\vspace*{0pt}\label{z5}%
\end{equation}and
	\begin{equation}
		\begin{array}
			[c]{cl}
			\text{{\small \text{$\mathcal{P}$}(\textit{n}) \textit{\ }= }} & \text{+\thinspace
				{\small 250765325100000\thinspace}}{\small n}^{11}\,\text{{\small --\thinspace
					1087318449630000\thinspace}}{\small n}^{10}\,\text{{\small +\thinspace
					2067749814046250\thinspace}}{\small n}^{9}%
			\genfrac{}{}{0pt}{0}{\genfrac{}{}{0pt}{1}{{}}{{}}}{{}}%
			\\
			& \text{{\small --\thinspace2269551612681475\thinspace}}{\small n}%
			^{8}\,\text{{\small +\thinspace1592180015776565\thinspace}}{\small n}%
			^{7}\,\text{{\small --\thinspace746938801646725\thinspace}}{\small n}^{6}\,%
			\genfrac{}{}{0pt}{0}{\genfrac{}{}{0pt}{1}{{}}{{}}}{{}}%
			\\
			& \text{{\small +\thinspace238210943593421\thinspace}}{\small n}%
			^{5}\,\text{{\small --\thinspace51452348050672\thinspace}}{\small n}%
			^{4}\,\text{{\small +\thinspace7352050259484\thinspace}}{\small n}^{3}\,%
			\genfrac{}{}{0pt}{0}{\genfrac{}{}{0pt}{1}{{}}{{}}}{{}}%
			\\
			& \text{{\small --\thinspace660416507568\thinspace}}{\small n}^{2}%
			\,\text{{\small +\thinspace33552610560\thinspace}}{\small n}%
			\,\,\text{{\small --\thinspace731566080}}%
			\genfrac{}{}{0pt}{0}{\genfrac{}{}{0pt}{1}{{}}{{}}}{{}}%
			\\
			\text{{\small \text{$\mathcal{R}$}(\textit{n}) \textit{\ }= }} &
			\text{{\small 24\thinspace\textit{n}}}^{5}\text{{\small (2\textit{n}--1)(3\textit{n}%
					--1)(3\textit{n}--2)(4\textit{n}--1)(4\textit{n}--3)(5\textit{n}%
					--1)(5\textit{n}--2)(5\textit{n}--3)(5\textit{n}--4)}}%
			\genfrac{}{}{0pt}{0}{\genfrac{}{}{0pt}{1}{{}}{{}}}{{}}%
		\end{array}
		\label{z6}%
	\end{equation}\vspace{0pt}

Series \textbf{\textit{i})} --Eqs.(\ref{z2}--\ref{z21}--\ref{z3})-- is a very fast auxiliary algorithm taking part to compute some Dirichlet \textit{L}--function formulas $L_k(s)$ for $s=3$ as it is shown further on. It has been incorporated in Number Theory software library FLINT \cite{FLINT0} as the standard algorithm to calculate constant $\zeta(3)$ (see documentation section 9.31.5). Series \textbf{\textit{i})} for computing and series \textbf{\textit{ii})} --Eqs.(\ref{z2}--\ref{z5}--\ref{z6})-- for digits validation, were applied to extend the number of known decimal places for this constant beyond 2\textit{e}+12 digits \cite{YEE0} (2023).  

\pagebreak
\subsubsection{$\boldsymbol{L_{-3}(2)}\ $constant}This is a self-dual primitive \textit{L}--function value that was recently proven to be irrational \cite{CALE}. This primary constant takes place in formulas for some special \textit{L}--Series, particularly $L_{k}(s)$ for $s=2$ and characters $k=-15,-24$. There are currently two known fast rational algorithms of hypergeometric--type to compute $L_{-3}(2)$.\vspace{8pt} 

The first one is a conjectured identity ($\bumpeq$) from Z. W. Sun (2023) \cite{SUN} (Conjecture 5.1)
\begin{equation}\label{zw}
	L_{-3}(2)\bumpeq\sum_{n=1}^\infty\,\frac{20\,n-6}{n^2\,(2\,n-1)}\cdot\left(\frac1{81}\right)^n\left[\begin{array}
		[l]{cc}%
		1 & 1\\
		\frac13 & \frac23 
	\end{array}\right]_n
\end{equation}This series has cost $C_s=2.7307...$\vspace{8pt} 

The second algorithm is Guillera's fast proven series (\cite{GUILLERA}, see examples in Maple\texttrademark$\ $code). It is obtained from Eq.(\ref{9y}) using  $L_{-3}(2)=\tfrac19\,\mathcal{S}_0(\tfrac13,\tfrac23)$. In this case some Pochhammer's simplify, deepness falls down to ${\small \textit{D}}=5$ and cost reduces to $C_s=20/\log(1024)=2.8854..$.\vspace{0pt} 
More than 1\textit{e}+12 decimal places are currently known for this constant (2025) \cite{YEE0}. Sun's series was applied in June 2025 for computing and Guillera's for digits validation.

\subsubsection{$\boldsymbol{L_{-4}(2)},$ Catalan's constant $\boldsymbol{G}$}This primary constant participates in $L_{-20}(2)$ formula. $G$ is usually computed via Gosper--Pilehrood series. This is an efficient hypergeometric algorithm found by Bill Gosper in the seventies for the MIT--MACSYMA project --\cite{BRAD} Eq.(77)--. Series was rediscovered more than 30 years later (2010) by the Pilehroods --\cite{HESS} Corollary 2--. It is found and proven applying Eq.(\ref{6}) with command 
\begin{equation}\label{g1}
	\texttt{SER(}1 + 4 n, \tfrac12 + 2 n, \tfrac12 + n, \tfrac12 +  n, \texttt{oo)} 
\end{equation}It is also obtained using $G=\tfrac1{16}\,\mathcal{\bar{S}}_0(\tfrac14,\tfrac34)$ where $\mathcal{\bar{S}}_0(x,y)$ is the difference of Hurwitz $\zeta$'s, Eq.(\ref{8xy}), accelerated by applying Dougall's search pattern $(\alpha,\beta,\gamma,\delta,\varepsilon)$ = $(4,1,1,1,1)$.
\begin{equation}\label{g10}
	\texttt{SER(}1 + 4 n, \tfrac14 + n, \tfrac14 + n, \tfrac34 + n, \tfrac34 + n \texttt{)} 
\end{equation}provides a very fast condensed series that sums two consecutive terms at once so that the original Gosper--Pilehrood series is retrieved by bisection. The binary splitting cost of this algorithm is $C_s=16/\log(729/4)=3.07374...$\vspace{4pt}

Another very efficient Catalan's constant series was also found via WZ--Dougall $_5H_5$ method by this author (2023). This is proven with command
\begin{equation}\label{g2}
	\texttt{SER(}1 + 8 n, \tfrac12 + 4 n, \tfrac12 + 3 n, \tfrac12 + 2 n, \texttt{oo)} 
\end{equation}The resulting series has $\rho=-1/12500$ giving more than 4 decimal places per term. It has been incorporated in Number Theory software library FLINT \cite{FLINT0} as the standard algorithm to compute constant $G$ (see documentation section 9.31.4). It is given by
\begin{equation}
	G=
	{\displaystyle\sum\limits_{n=1}^{\infty}}\ \left(-\frac1{2^2\,5^5}\right)^n\cdot
	\frac{\mathcal{P}(n)}{\mathcal{R}(n)}\cdot\left[\begin{array}
		[l]{cccccccc}%
		1 & 1 & 1 & \frac12 & \frac13 & \frac23 & \frac16 & \frac56 
		\genfrac{}{}{0pt}{1}{\genfrac{}{}{0pt}{0}{{}}{{}}}{{}}%
		\\
		\frac1{10} & \frac3{10} & \frac{7}{10} & \frac{9}{10} & \frac{1}{12} &
		\frac{5}{12} & \frac{7}{12} & \frac{11}{12} 
		\genfrac{}{}{0pt}{1}{\genfrac{}{}{0pt}{0}{{}}{{}}}{{}}%
	\end{array}\right]_n\vspace*{-8pt}\label{g3}%
\end{equation}
\begin{equation}
	\begin{array}
		[c]{cl}
		& \\
		\text{{\small \text{$\mathcal{P}$}(\textit{n})\textit{\ }=}} &
		\,\text{{\small +\thinspace43203456\thinspace\textit{n}}}^{6}%
		\,\text{{\small --\thinspace92809152\thinspace\textit{n}}}^{5}%
		\,\text{{\small +\thinspace76613904\thinspace\textit{n}}}^{4}
		\,\text{{\small --\thinspace30494304\thinspace\textit{n}}}^{3}
		\genfrac{}{}{0pt}{0}{\genfrac{}{}{0pt}{1}{{}}{{}}}{{}}%
		\\
		& \,\text{{\small +\thinspace6004944\thinspace\textit{n}}}^{2}%
		\,\text{{\small --\thinspace536620\thinspace\textit{n}}}%
		\,\text{{\small +\thinspace17325}}%
		\genfrac{}{}{0pt}{0}{\genfrac{}{}{0pt}{1}{{}}{{}}}{{}}%
		\\
		\text{{\small \text{$\mathcal{R}$}(\textit{n}) \textit{\ }= }} &
		\text{{\small --\thinspace768\thinspace\textit{n}}}^{3}\text{{\small (2\textit{n}--1)(3\textit{n}%
				--1)(3\textit{n}--2)(6\textit{n}--1)(6\textit{n}--5)}}%
		\genfrac{}{}{0pt}{0}{\genfrac{}{}{0pt}{1}{{}}{{}}}{{}}%
	\end{array}
	\label{g4}%
\end{equation}\vspace{0pt}Observe that $\pi,\ \zeta(3),\ L_{-3}(2),\ G$ are sorted from fast to slower computing constants, so the $L_{k}(s)$ formulas that only require the faster ones are certainly more efficient. In the next section these fast Dirichlet's $L_{k}(s)$ formulas are introduced.       

\subsection{Fast $\boldsymbol{L_{k}(s)}$ Formulas}\hspace{-3pt}Eqs.(\ref{8s}--\ref{8w}--\ref{8t}--\ref{8u}--\ref{8v}) plus Eq.(\ref{L2}) together with $\ \mathcal{S}_0(x,y)$ and $\mathcal{S}_1(x)$ provide all the pieces to lift up new proven identities for these selected primitive Dirichlet $L$--Series: $L_{-7}(2)$, $L_{-8}(2)$, $L_{-15}(2)$, $L_{-20}(2)$,  $L_{-24}(2)$, $L_5(3)$, $L_8(3)$ and $L_{12}(3)$.\vspace{8pt} 

Note that for some particular $x,y$ values, $\mathcal{S}_0(x,y)$ --Eq.(\ref{9y})-- or $\mathcal{S}_1(x)$ --Eq.(\ref{10y})-- might be replaced by even faster series $\bar{\mathcal{S}}_0(x,y)$ or $\bar{\mathcal{S}}_1(x)$ if other search patterns $(\alpha,\beta,\gamma,\delta,\varepsilon)$ give higher accelerations. This happens for $L_{-7}(2)$ whose best pattern is $(5,2,2,2,2)$, $L_{-8}(2)$, $(4,1,1,1,1)$ and $L_{5}(3)$ with $(6,3,2,2,2)$ as it is seen in what follows.
\subsubsection{$\boldsymbol{L_{-7}(2)}$}Apparently there is currently just one known linearly convergent rational algorithm of hypergeometric–type to compute this $L$--function value. Guillera's conjectured series \cite{GUI2} Eq.(28) derived as the dual of a divergent hypergeometric series for $\pi^{-1}$,\vspace{-2pt} 
\begin{equation}\label{L72}
	L_{-7}(2)\bumpeq\frac1{7^4}\cdot\sum_{n=1}^\infty\,\frac{2210\,n^2-1273\,n+120}{-n^3\,(2\,n-1)}\cdot\left(-\frac{7^4}{2^{14}}\right)^n\left[\begin{array}
		[l]{cccc}%
		\tfrac12 & 1 & 1 & 1\vspace{2pt}\\
		\frac18 & \frac38 & \frac58 & \frac78\vspace{2pt} 
	\end{array}\right]_n
\end{equation}this dual series has binary splitting cost $C_s=16/\log(16384/2401)=8.3315...$\vspace{8pt}\newline
A different identity is obtained from Eq.(\ref{L2}) where the $L$--function value $L_{-7}(2)$ is expressed as a linear combination of three linearly convergent accelerated series $\bar{\mathcal{S}}_0(x,y),\ $\vspace{0pt}
\begin{equation}\label{L72a}
	49\,L_{-7}(2)=\,\bar{\mathcal{S}}_0(\tfrac17,\tfrac67)+\bar{\mathcal{S}}_0(\tfrac27,\tfrac57)-\bar{\mathcal{S}}_0(\tfrac37,\tfrac47)\vspace{0pt}	
\end{equation}By using the search pattern $(\alpha,\beta,\gamma,\delta,\varepsilon)=(5,2,2,2,2)$\vspace{2pt} and applying these commands
\begin{equation}\label{L72c}
	\texttt{SER(}1 + 5 n, \tfrac\ell7 + 2 n, \tfrac\ell7 + 2 n, \tfrac{7-\ell}7 + 2 n,  \tfrac{7-\ell}7 + 2 n\texttt{)}\ \ \ \ell=1,2,3\vspace{1pt} 
\end{equation}these series are written employing Eq.(\ref{7a}) convention with $\rho=-64/531441$
\begin{equation}\label{L72e}
	\bar{\mathcal{S}}_0(\tfrac{\ell}{7},\tfrac{7-\ell}{7}) = \sum_{n=1}^\infty\,\frac{\mathcal{P}_\ell(n)}{\mathcal{R}_\ell(n)}\cdot\left(-\frac{2^6}{3^{12}}\right)^n\cdot\mathcal{M}_n^{(\ell)},\ \ \ \ell=1,2,3 
\end{equation}
\begin{equation*}
\begin{array}
		[l]{ccl}
		\text{{\small \negthinspace\hspace{-18pt}\text{$\mathcal{P}_1$}}}\text{{\small (\textit{n}) \textit{\ \ } }} & \space\space= &\text{+\thinspace
			{\small 3299714444390\thinspace}}{\small n}^{9}\,\text{{\small --\thinspace
				11086179271881\thinspace}}{\small n}^{8}\,\text{{\small +\thinspace
				15792915284685\thinspace}}{\small n}^{7}%
		\genfrac{}{}{0pt}{0}{\genfrac{}{}{0pt}{1}{{}}{{}}}{{}}%
		\\
		& &\text{{\small --\thinspace12436203664563\thinspace}}{\small n}%
		^{6}\,\text{{\small +\thinspace5906387588241\thinspace}}{\small n}%
		^{5}\,\text{{\small --\thinspace1727982499236\thinspace}}{\small n}^{4}\,%
		\genfrac{}{}{0pt}{0}{\genfrac{}{}{0pt}{1}{{}}{{}}}{{}}%
		\\
		& &\text{{\small +\thinspace303911726860\thinspace}}{\small n}%
		^{3}\,\text{{\small --\thinspace29742768720\thinspace}}{\small n}%
		^{2}\,\text{{\small +\thinspace1376968320\thinspace}}{\small n}\,\text{{\small --\thinspace23362560\thinspace}}\vspace{2pt}%
		\genfrac{}{}{0pt}{0}{\genfrac{}{}{0pt}{1}{{}}{{}}}{{}}%
		\\
		\text{{\small \negthinspace\hspace{-18pt}\text{$\mathcal{R}_1$}}}\text{{\small (\textit{n}) \textit{\ \ }}} & \space\space= &
		\text{{\small --\thinspace56\thinspace}{\small \textit{n}}}^{3}\text{{\small \thinspace[\thinspace(7\textit{n}--4)(7\textit{n}%
				--3)(14\textit{n}--13)(14\textit{n}--1)\thinspace]}}^2%
	\end{array}
	\end{equation*}
\begin{equation}\label{L72f}%
	\begin{array}
		[l]{ccl}
		\text{{\small \space\text{$\mathcal{P}_2$}}}\text{{\small (\textit{n}) \textit{\ } }} & = &\text{+\thinspace
			{\small 3299714444390\thinspace}}{\small n}^{9}\,\text{{\small --\thinspace
				12028954827421\thinspace}}{\small n}^{8}\,\text{{\small +\thinspace
				18694624573863\thinspace}}{\small n}^{7}%
		\genfrac{}{}{0pt}{0}{\genfrac{}{}{0pt}{1}{{}}{{}}}{{}}%
		\\
		& &\text{{\small --\thinspace16165440254775\thinspace}}{\small n}%
		^{6}\,\text{{\small +\thinspace8510524572599\thinspace}}{\small n}%
		^{5}\,\text{{\small --\thinspace2805337398948\thinspace}}{\small n}^{4}\,%
		\genfrac{}{}{0pt}{0}{\genfrac{}{}{0pt}{1}{{}}{{}}}{{}}%
		\\
		& &\text{{\small +\thinspace573713306964\thinspace}}{\small n}%
		^{3}\,\text{{\small --\thinspace69548120208\thinspace}}{\small n}%
		^{2}\,\text{{\small +\thinspace4497806880\thinspace}}{\small n}\,\text{{\small --\thinspace118195200\thinspace}}\vspace{2pt}%
		\genfrac{}{}{0pt}{0}{\genfrac{}{}{0pt}{1}{{}}{{}}}{{}}%
		\\
		\text{{\small \space\text{$\mathcal{R}_2$}}}\text{{\small (\textit{n}) \textit{\ } }} & = &
		\text{{\small --\thinspace56\thinspace}{\small \textit{n}}}^{3}\text{{\small \thinspace[\thinspace(7\textit{n}--6)(7\textit{n}%
				--1)(14\textit{n}--9)(14\textit{n}--5)\thinspace]}}^2%
	\end{array}	
\end{equation}
\begin{equation*} 
	\begin{array}
		[l]{ccl}
		\text{{\small \hspace{-1pt}\text{$\mathcal{P}_3$}}}\text{{\small (\textit{n}) \textit{\ } }}& = & \text{+\thinspace
			{\small 3299714444390\thinspace}}{\small n}^{9}\,\text{{\small --\thinspace
				12971730382961\thinspace}}{\small n}^{8}\,\text{{\small +\thinspace
				22065924669985\thinspace}}{\small n}^{7}%
		\genfrac{}{}{0pt}{0}{\genfrac{}{}{0pt}{1}{{}}{{}}}{{}}%
		\\
		& &\text{{\small --\thinspace21255429649979\thinspace}}{\small n}%
		^{6}\,\text{{\small +\thinspace12728430435805\thinspace}}{\small n}%
		^{5}\,\text{{\small --\thinspace4891758484592\thinspace}}{\small n}^{4}\,%
		\genfrac{}{}{0pt}{0}{\genfrac{}{}{0pt}{1}{{}}{{}}}{{}}%
		\\
		& &\text{{\small +\thinspace1200733928744\thinspace}}{\small n}%
		^{3}\,\text{{\small --\thinspace180736839864\thinspace}}{\small n}%
		^{2}\,\text{{\small +\thinspace15099749280\thinspace}}{\small n}\,\text{{\small --\thinspace533174400\thinspace}}\vspace{2pt}%
		\genfrac{}{}{0pt}{0}{\genfrac{}{}{0pt}{1}{{}}{{}}}{{}}%
		\\
		\text{{\small \hspace{-1pt}\text{$\mathcal{R}_3$}}}\text{{\small (\textit{n}) \textit{\ }}} & = &
		\text{{\small --\thinspace56\thinspace}{\small \textit{n}}}^{3}\text{{\small \thinspace[\thinspace(7\textit{n}--5)(7\textit{n}%
				--2)(14\textit{n}--11)(14\textit{n}--3)\thinspace]}}^2%
	\end{array}%
\end{equation*}and
\begin{equation*}
	\mathcal{M}_n^{(1)}=\left[\begin{array}
		[l]{ccccccccccc}%
		1 & 1 & 1 & \frac37  & \frac37  & \frac47 & \frac47 & \frac1{14} & \frac1{14} & \frac{13}{14} & \frac{13}{14}   
		\genfrac{}{}{0pt}{1}{\genfrac{}{}{0pt}{0}{{}}{{}}}{{}}%
		\\
		\frac12 & \frac27 & \frac57 & \frac1{21} & \frac1{21}  & \frac8{21}  & \frac8{21}  & \frac{13}{21} & \frac{13}{21} & \frac{20}{21} & \frac{20}{21} 
		\genfrac{}{}{0pt}{1}{\genfrac{}{}{0pt}{0}{{}}{{}}}{{}}%
	\end{array}\right]_n
\end{equation*}
\begin{equation}\label{L72h}
	\mathcal{M}_n^{(2)}=\left[\begin{array}
		[l]{ccccccccccc}%
		1 & 1 & 1 & \frac17  & \frac17  & \frac67 & \frac67 & \frac5{14} & \frac5{14} & \frac{9}{14} & \frac{9}{14}   
		\genfrac{}{}{0pt}{1}{\genfrac{}{}{0pt}{0}{{}}{{}}}{{}}%
		\\
		\frac12 & \frac37 & \frac47 & \frac2{21} & \frac2{21}  & \frac5{21}  & \frac5{21}  & \frac{16}{21} & \frac{16}{21} & \frac{19}{21} & \frac{19}{21} 
		\genfrac{}{}{0pt}{1}{\genfrac{}{}{0pt}{0}{{}}{{}}}{{}}%
	\end{array}\right]_n
\end{equation}
\begin{equation*}
	\mathcal{M}_n^{(3)}=\left[\begin{array}
		[l]{ccccccccccc}%
		1 & 1 & 1 & \frac27  & \frac27  & \frac57 & \frac57 & \frac3{14} & \frac3{14} & \frac{11}{14} & \frac{11}{14}   
		\genfrac{}{}{0pt}{1}{\genfrac{}{}{0pt}{0}{{}}{{}}}{{}}%
		\\
		\frac12 & \frac17 & \frac67 & \frac4{21} & \frac4{21}  & \frac{10}{21}  & \frac{10}{21}  & \frac{11}{21} & \frac{11}{21} & \frac{17}{21} & \frac{17}{21} 
		\genfrac{}{}{0pt}{1}{\genfrac{}{}{0pt}{0}{{}}{{}}}{{}}%
	\end{array}\right]_n\vspace{8pt}
\end{equation*}The binary splitting cost for this formula is $C_s=3\cdot44/\log(531441/64)=14.6269...$, too large compared to Guillera's series Eq.(\ref{L72}), however this new proven tri--series identity is useful anyway because it is applied to validate the large number of decimal places that is obtained by the fast conjectural series.      
Until before this note was written, there was not a fast secondary identity for this $L$--Series. Eq.(\ref{L72a}) is the unique irreducible formula requiring multiple series in this work.
\subsubsection{$\boldsymbol{L_{-8}(2)}$}For this $L$--Series, the expressions
\begin{equation}\label{L82a}
32\,L_{-8}(2)=\zeta(2,\tfrac38)-\zeta(2,\tfrac78)+2\,\pi^2\sqrt2	
\end{equation}
\begin{equation}\label{L82b}
32\,L_{-8}(2)=\zeta(2,\tfrac18)-\zeta(2,\tfrac58)-2\,\pi^2\sqrt2\vspace{4pt}
\end{equation}give two algorithmically distinct identities where each difference of $\zeta$'s for $(x,y)=(\tfrac38,\tfrac78)$ and $(x,y)=(\tfrac18,\tfrac58)$ is accelerated either by $\mathcal{S}_0(x,y)$ --Eq.(\ref{9y})-- with cost $C_s=4.6166...$ or by using the search pattern  
($\alpha,\beta,\gamma,\delta,\varepsilon$)=(4,1,1,1,1) that brings $\bar{\mathcal{S}}_0(x,y)$ with cost $C_s=48/\log(531441/16)=4.6106...$ a slightly faster series. In this case
\begin{equation}\label{L82c}
	32\,L_{-8}(2)=\bar{\mathcal{S}}_0(\tfrac38,\tfrac78)+2\,\pi^2\sqrt2	
\end{equation}
\begin{equation}\label{L82d}
	32\,L_{-8}(2)=\bar{\mathcal{S}}_0(\tfrac18,\tfrac58)-2\,\pi^2\sqrt2\vspace{-4pt}
\end{equation}with    
\begin{equation}\label{L82e}
	\bar{\mathcal{S}}_0(x,y)=\sum_{n=1}^\infty\,\frac{\mathcal{P}(n,x,y)}{\mathcal{R}(n,x,y)}\cdot\rho^n\cdot\mathcal{M}_n(x,y)\vspace{0pt}
\end{equation}where $\mathcal{P},\mathcal{R},\mathcal{M}_n$ and $\rho=2^4/3^{12}$ are obtained with the corresponding commands\vspace{4pt}
\begin{equation}\label{L82f}\vspace{1pt}
	\texttt{SER(}\tfrac54 + 4\,n, \tfrac38 + n, \tfrac38 + n, \tfrac78 + n,  \tfrac78 + n\texttt{)} 
\end{equation}
\begin{equation}\label{L82g}\vspace{1pt}
	\texttt{SER(}\tfrac34 + 4\,n, \tfrac18 + n, \tfrac18 + n, \tfrac58 + n,  \tfrac58 + n\texttt{)} 
\end{equation}it gives\vspace{-2pt} 
\begin{equation}\label{L82h}
	\begin{array}
		[l]{ccl}
		\text{{\small \negthinspace\hspace{-20pt}\text{$\mathcal{P}(\text{\textit{n}},\tfrac38,\tfrac78)$}}} & \space\space= &\text{+\thinspace
			{\small 114122649763840\thinspace}}{\small n}^{10}\,\text{{\small --\thinspace
				429099445125120\thinspace}}{\small n}^{9}\,\text{{\small +\thinspace
				693932983844864\thinspace}}{\small n}^{8}%
		\genfrac{}{}{0pt}{0}{\genfrac{}{}{0pt}{1}{{}}{{}}}{{}}%
		\\
		& &\text{{\small --\thinspace631504224911360\thinspace}}{\small n}^{7}\,
		\text{{\small +\thinspace355139912728576\thinspace}}{\small n}^{6}\,
		\text{{\small --\thinspace127514032996352\thinspace}}{\small n}^{5}\,%
		\genfrac{}{}{0pt}{0}{\genfrac{}{}{0pt}{1}{{}}{{}}}{{}}%
		\\
		& &\text{{\small +\thinspace29138090287104\thinspace}}{\small n}^{4}\,
		\text{{\small --\thinspace4087411945472\thinspace}}{\small n}^{3}\,
		\text{{\small +\thinspace324916373312\thinspace}}{\small n}^{2}\,
		\genfrac{}{}{0pt}{0}{\genfrac{}{}{0pt}{1}{{}}{{}}}{{}}
		\\
		& & \text{{\small --\thinspace12542210880\thinspace}}{\small n}\,
		\text{{\small +\thinspace180650925\thinspace}}\,
		\genfrac{}{}{0pt}{0}{\genfrac{}{}{0pt}{1}{{}}{{}}}{{}}
		\\
		\text{{\small \negthinspace\hspace{-20pt}\text{$\mathcal{R}(\text{\textit{n}},\tfrac38,\tfrac78)$}}} & \space\space= &
		\text{{\small \thinspace65536\thinspace}{\small [\thinspace\textit{n}\thinspace(2\textit{n}--1)\thinspace]}}^{3}
		\text{{\small \thinspace(4\textit{n}--1)(4\textit{n}--3)\thinspace[\thinspace(8\textit{n}--1)(8\textit{n}--5)\thinspace]}}^{2}%
			%
	\end{array}
\end{equation}
\begin{equation}\label{L82i}
	\begin{array}
		[l]{ccl}
		\text{{\small \negthinspace\hspace{-20pt}\text{$\mathcal{P}(\text{\textit{n}},\tfrac18,\tfrac58)$}}} & \space\space= &\text{+\thinspace
			{\small 114122649763840\thinspace}}{\small n}^{10}\,\text{{\small --\thinspace
				520397564936192\thinspace}}{\small n}^{9}\,\text{{\small +\thinspace
				1036756736540672\thinspace}}{\small n}^{8}%
		\genfrac{}{}{0pt}{0}{\genfrac{}{}{0pt}{1}{{}}{{}}}{{}}%
		\\
		& &\text{{\small --\thinspace1185072896016384\thinspace}}{\small n}^{7}\,
		\text{{\small +\thinspace857946029817856\thinspace}}{\small n}^{6}\,
		\text{{\small --\thinspace409505218428928\thinspace}}{\small n}^{5}\,%
		\genfrac{}{}{0pt}{0}{\genfrac{}{}{0pt}{1}{{}}{{}}}{{}}%
		\\
		& &\text{{\small +\thinspace129930258636800\thinspace}}{\small n}^{4}\,
		\text{{\small --\thinspace26919161769984\thinspace}}{\small n}^{3}\,
		\text{{\small +\thinspace3464815945536\thinspace}}{\small n}^{2}\,
		\genfrac{}{}{0pt}{0}{\genfrac{}{}{0pt}{1}{{}}{{}}}{{}}
		\\
		& & \text{{\small --\thinspace248697546048\thinspace}}{\small n}\,
		\text{{\small +\thinspace7536198285\thinspace}}\,\vspace{2pt}
		\genfrac{}{}{0pt}{0}{\genfrac{}{}{0pt}{1}{{}}{{}}}{{}}
		\\
		\text{{\small \negthinspace\hspace{-20pt}\text{$\mathcal{R}(\text{\textit{n}},\tfrac18,\tfrac58)$}}} & \space\space= &
		\text{{\small \thinspace65536\thinspace}{\small [\thinspace\textit{n}\thinspace(2\textit{n}--1)\thinspace]}}^{3}
		\text{{\small \thinspace(4\textit{n}--1)(4\textit{n}--3)\thinspace[\thinspace(8\textit{n}--3)(8\textit{n}--7)\thinspace]}}^{2}\vspace{2pt}
		\vspace{0pt}
	\end{array}
\end{equation}and\vspace{-2pt}
\begin{equation}\label{L82j}
	\mathcal{M}_n(\tfrac38,\tfrac78)=\left[\begin{array}
		[l]{cccccccccccc}%
		1 & 1 & 1 & \frac12  & \frac12  & \frac12 & \frac14 & \frac34 & \frac38 & \frac38 & \frac78 & \frac78  
		\genfrac{}{}{0pt}{1}{\genfrac{}{}{0pt}{0}{{}}{{}}}{{}}%
		\\
		\frac18 & \frac18 & \frac58 & \frac58 & \frac7{24}  & \frac{7}{24}  & \frac{11}{24}  & \frac{11}{24} & \frac{19}{24} & \frac{19}{24} & \frac{23}{24} & \frac{23}{24} 
		\genfrac{}{}{0pt}{1}{\genfrac{}{}{0pt}{0}{{}}{{}}}{{}}%
	\end{array}\right]_n
\end{equation}\vspace{0pt}
\begin{equation}\label{L82k}
	\mathcal{M}_n(\tfrac18,\tfrac58)=\left[\begin{array}
		[l]{cccccccccccc}%
		1 & 1 & 1 & \frac12  & \frac12  & \frac12 & \frac14 & \frac34 & \frac18 & \frac18 & \frac58 & \frac58  
		\genfrac{}{}{0pt}{1}{\genfrac{}{}{0pt}{0}{{}}{{}}}{{}}%
		\\
		\frac38 & \frac38 & \frac78 & \frac78 & \frac1{24}  & \frac{1}{24}  & \frac{5}{24}  & \frac{5}{24} & \frac{13}{24} & \frac{13}{24} & \frac{17}{24} & \frac{17}{24} 
		\genfrac{}{}{0pt}{1}{\genfrac{}{}{0pt}{0}{{}}{{}}}{{}}%
	\end{array}\right]_n\vspace{4pt}
\end{equation}By using these series each identity in Eqs.(\ref{L82c}--\ref{L82d}) has a total binary splitting cost $C_s=48/\log(531441/16)$ ($\bar{\mathcal{S}}_0$ series) $+\ 0.36748$ ($\pi$ Chudnovsky) = 4.9781... giving the most efficient new formulas in this work. 
\subsubsection{$\boldsymbol{L_{-15}(2)}$}This $L$--Series satisfies the following computationally independent identities that are proven from Eq.(\ref{L2}), applying the reduction formulas, taking $x=\tfrac1{15}$ and $k=5$ in Eq.(\ref{8u}) and replacing the difference of $\zeta$'s by the fast series $\mathcal{S}_0$ from Eq.(\ref{9y}),
\begin{equation}\label{L15a}
	225\,L_{-15}(2)=4\,\mathcal{S}_0(\tfrac{1}{15},\tfrac{11}{15})-234\thinspace L_{-3}(2)-8\thinspace\pi^2\sqrt{15+6\sqrt5}\vspace{-4pt}	
\end{equation}
\begin{equation}\label{L15c}
	225\,L_{-15}(2)=4\,\mathcal{S}_0(\tfrac{2}{15},\tfrac{7}{15})+234\thinspace L_{-3}(2)-8\thinspace\pi^2\sqrt{15-6\sqrt5}\vspace{4pt}	
\end{equation}Each identity is computed with a total binary splitting cost $C_s=5.19372...$ (series $\mathcal{S}_0$) + $2.73072...$ ($L_{-3}(2)$, Z.W. Sun) + $0.36748...$ ($\pi$, Chudnovsky) = 8.2919...\vspace{4pt}\newline Note that by mirroring $\mathcal{S}_0(x,y)$ with $(x,y)\in\{(\tfrac{4}{15},\tfrac{14}{15}),(\tfrac{8}{15},\tfrac{13}{15})\}$, Eq.(\ref{8w}) together with Eqs.(\ref{L15a}--\ref{L15c}) produce two additional similar identities employing different series\footnote{In this case the mirror identities are $225\,L_{-15}(2)=4\,\mathcal{S}_0(\tfrac{4}{15},\tfrac{14}{15})-234\thinspace L_{-3}(2)+8\thinspace\pi^2\sqrt{15+6\sqrt5}$ and $225\,L_{-15}(2)=4\,\mathcal{S}_0(\tfrac{8}{15},\tfrac{13}{15})+234\thinspace L_{-3}(2)+8\thinspace\pi^2\sqrt{15-6\sqrt5}$. See \textit{Appendix} for mirror $L$--Series identities}.
\subsubsection{$\boldsymbol{L_{-20}(2)}$}The next algorithmically different identities are derived from Eq.(\ref{L2}) applying Hurwitz $\zeta$'s reduction formulas, by taking $x=\tfrac1{20}$ and $k=5$ in Eq.(\ref{8u}) and series $\mathcal{S}_0$ from Eq.(\ref{9y}). $G$ is Catalan's constant.\vspace{-8pt}
\begin{equation}\label{L20a}
	100\,L_{-20}(2)=\mathcal{S}_0(\tfrac{1}{20},\tfrac{11}{20})-96\thinspace G-2\thinspace\pi^2\sqrt{50+22\sqrt5}\vspace{-2pt}	
\end{equation}
\begin{equation}\label{L20c}
	100\,L_{-20}(2)=\mathcal{S}_0(\tfrac{3}{20},\tfrac{13}{20})+96\thinspace G-2\thinspace\pi^2\sqrt{50-22\sqrt5}\vspace{8pt}
\end{equation}Each identity $L_{-20}(2)$ has a total binary splitting cost $C_s=4.61662...$ (series $\mathcal{S}_0$) + $3.07374...$ ($G$, Gosper--Pilehrood) + $0.36748...$ ($\pi$, Chudnovsky) = 8.0578...\vspace{6pt}\newline By taking $(x,y)=(\tfrac{9}{20},\tfrac{19}{20})$ and $(x,y)=(\tfrac{7}{20},\tfrac{17}{20})$, Eq.(\ref{8w}) together with Eqs.(\ref{L20a}--\ref{L20c}) give two additional similar mirror identities with different series \vspace{8pt}$\mathcal{S}_0(x,y)$. (See \textit{Appendix}).

\subsubsection{$\boldsymbol{L_{-24}(2)}$}Using Eq.(\ref{L2}) and applying $\mathcal{S}_0$ from Eq.(\ref{9y}) this $L$--Series has two kind of computationally different identities depending on the ancillary constant used. By taking $x=\tfrac1{24}$ and $k=8$ in Eq.(\ref{8u}) it gives $L_{-3}(2)$ and
\begin{equation}\label{L24a}
	144\,L_{-24}(2)=\mathcal{S}_0(\tfrac{1}{24},\tfrac{17}{24})-180\thinspace L_{-3}(2)-2\thinspace\pi^2(4\sqrt3+3\sqrt6)	
\end{equation}
\begin{equation}\label{L24b}
	144\,L_{-24}(2)=\mathcal{S}_0(\tfrac{5}{24},\tfrac{13}{24})+180\thinspace L_{-3}(2)+2\thinspace\pi^2(4\sqrt3-3\sqrt6)\vspace{4pt}
\end{equation}Binary splitting cost for these $L_{-24}(2)$ identities is the same as $L_{-15}(2)$, $C_s=8.2919...$ In this case mirror identities are derived from Eq.(\ref{8w}) by using $(x,y)\in\{(\tfrac{7}{24},\tfrac{23}{24}),(\tfrac{11}{24},\tfrac{19}{24})\}$.
If $x=\tfrac1{24}$ and $k=6$ is taken in Eq.(\ref{8u}) Catalan's constant $G$ rises as a linear combination of Hurwitz $\zeta$'s that are reduced producing these algorithmically independent identities
\begin{equation}\label{L24c}
	144\,L_{-24}(2)=\mathcal{S}_0(\tfrac{5}{24},\tfrac{23}{24})-160\thinspace G+2\thinspace\pi^2(5\sqrt2+4\sqrt3)\vspace{4pt}	
\end{equation}
\begin{equation}\label{L24d}
	144\,L_{-24}(2)=\mathcal{S}_0(\tfrac{11}{24},\tfrac{17}{24})+160\thinspace G+2\thinspace\pi^2(5\sqrt2-4\sqrt3)\vspace{4pt}
\end{equation}with a higher total cost $C_s=8.6349...$ Note that by taking $(x,y)\in\{(\tfrac{1}{24},\tfrac{19}{24}),(\tfrac{7}{24},\tfrac{13}{24})\}$ two additional similar mirror identities employing different series are obtained.

\subsubsection{$\boldsymbol{L_{5}(3)}$}For this special $L$--function value an accelerated series $\bar{\mathcal{S}}_1(x)=\zeta(3,x)$ is proven by applying WZ--Dougall's $_5H_5$ search pattern $(\alpha,\beta,\gamma,\delta,\varepsilon) = (6,3,2,2,2)$ with seed Eq.(\ref{8x}). The following algorithmically independent identities are obtained from Eq.(\ref{L2}) by using the reduction formulas Eq.(\ref{8t}) and Eq.(\ref{8v})
\begin{equation}\label{L53a}
	\negthinspace 625\,L_5(3)=\ 20\,\bar{\mathcal{S}}_1(\tfrac15)-620\,\zeta(3)-4\,\pi^3\left(\sqrt{25+2\sqrt5}+\sqrt{25-2\sqrt5}\right)\vspace{-6pt}
\end{equation}%
\begin{equation}\label{L53b}
625\,L_5(3)=-20\,\bar{\mathcal{S}}_1(\tfrac25)+620\,\zeta(3)+4\,\pi^3\left(\sqrt{25+2\sqrt5}-\sqrt{25-2\sqrt5}\right)
\end{equation}The accelerated Hurwitz $\zeta$ series are found with $x\in\{\tfrac15,\tfrac25\}$ (for additional accelerated mirror identities use Eq.(\ref{8t}) and take $x\in\{\tfrac35,\tfrac45\}$) applying the command\vspace{4pt} 
\begin{equation}\label{L53e}
\texttt{SER(}2\,x + 6\,n, x + 3 n, x + 2 n, x + 2 n, x + 2 n\texttt{)}\vspace{4pt} 	
\end{equation}This gives $\rho=-1/110592\ $ and
\begin{equation}\label{L53f}
	\bar{\mathcal{S}}_1(x)=\sum_{n=0}^\infty\,\frac{\mathcal{P}(n,x)}{\mathcal{Q}(n,x)}\cdot\left(-\frac1{2^{12}\,3^3}\right)^n\frac{(1)_n^5\,\left[\,(\frac12)_n(\frac{x}2)_n(\frac{x+1}2)_n\,\right]^3}{(\frac13)_n(\frac23)_n\left[\,(\frac{x+1}4)_n(\frac{x+2}4)_n(\frac{x+3}4)_n(\frac{x+4}4)_n\,\right]^3}
	\end{equation} 
\begin{flalign}\label{L53h}
	&\ \ \ \mathcal{Q}(n,x)=6\,x^3\,(3\,n + 2)(3\,n + 1)\left[(4\,n + x + 1)(4\,n + x + 2)(4\,n + x + 3)\,\right]^3\vspace{-8pt}&
\end{flalign}
\begin{multline*}
	\mathcal{P}(n,x) =\text{\small 16178176}\,n^{11}+\left(\text{\small 38134272}\,x \text{\small + 87494144}\right)\,n^{10}\\+\left(\text{\small 40486656}\, x^{2}\text{\small + 187222272}\, x \text{\small + 211522816}\right)\,n^{9}\\+\left(\text{\small 25499392}\, x^{3}\text{\small +178811136}\, x^{2}\text{\small +406560000}\, x \text{\small +301596160}\right)\, n^{8}\\+\left(\text{\small 10551552}\, x^{4}\text{\small +100197248}\, x^{3}\text{\small + 344790144}\, x^{2}\text{\small + 513893376}\, x \text{\small + 281718080}\right)\,n^{7}\\+\left(\text{\small 2997504}\, x^{5}\text{\small + 36382272}\,x^{4}\text{\small + 169112512}\, x^{3}\text{\small + 380623296}\, x^{2}\text{\small + 418517856}\, x\text{\small + 181010272}\right)\,n^{6}\\+(\text{\small 592048}\, x^{6}\text{\small + 8909424}\, x^{5}\text{\small + 52755504}\, x^{4}\text{\small + 159936816}\, x^{3}\\\text{\small + 264899280}\, x^{2}\text{\small + 229434480}\, x\text{\small + 81658480})\,n^{5}
\end{multline*}(cont.)\vspace{-6pt}
\begin{multline*}
+\vspace{-3pt}(\text{\small 80280}\, x^{7}\text{\small + 1481192}\, x^{6}\text{\small + 10820808}\, x^{5}\text{\small + 41637720}\, x^{4}\text{\small + 92597184}\, x^{3}\text{\small + 120473808}\, x^{2}\\ \text{\small +\vspace{-3pt} 85772400}\, x\text{\small + 25881872})\,n^{4}\\+\vspace{-3pt}(\text{\small 7152}\, x^{8}\text{\small + 163476}\, x^{7}\text{\small + 1452508}\, x^{6}\text{\small + 6860772}\, x^{5}\text{\small + 19288572}\, x^{4}\text{\small + 33574632}\, x^{3}\text{\small + 35803128}\, x^{2}\\ \text{\small +\vspace{-3pt} 21613440}\, x\text{\small + 5653120})\, n^{3}\\+\vspace{-3pt}(\text{\small 378}\, x^{9}\text{\small + 11268}\, x^{8}\text{\small + 122148}\, x^{7}\text{\small + 696278}\, x^{6}\text{\small + 2390292}\, x^{5}\text{\small + 5236914}\, x^{4}\text{\small + 7441070}\, x^{3}\\ \text{\small +\vspace{-3pt} 6709308}\, x^{2} \text{\small + 3519816}\, x\text{\small + 810864})\, n^{2}\\+\vspace{-1pt}(\text{\small 9}\, x^{10}\text{\small + 423}\, x^{9}\text{\small + 5769}\, x^{8}\text{\small + 39582}\, x^{7}\text{\small + 162757}\, x^{6}\text{\small + 433017}\, x^{5}\text{\small + 770661}\, x^{4}\text{\small + 921658}\, x^{3}\\ \text{\small +\vspace{-3pt} 720660}\, x^{2}\text{\small + 335448}\, x\text{\small + 68688})\, n \\ \text{\small +\vspace{-3pt} 6}\, x^{10}\text{\small + 114}\, x^{9}\text{\small + 957}\, x^{8}\text{\small + 4680}\, x^{7}\text{\small + 14804}\, x^{6}\text{\small + 31806}\, x^{5}\text{\small + 47385}\, x^{4}\\ \text{\small +\vspace{-3pt} 48896}\,x^{3}\text{\small + 33912}\, x^{2}\text{\small + 14256}\, x\text{\small + 2592}
\end{multline*}expressions valid for $x\in\{\tfrac15,\tfrac25,\tfrac35,\tfrac45\}$. The total binary splitting cost is $C_s=4.82193...$($\bar{\mathcal{S}}_1$ series) + $2.05136...$ ($\zeta(3)$, Eq.(\ref{z2})\thinspace + $0.36748...$ ($\pi$, Chudnovsky) = 7.2408... Certainly, it is also possible to use $\mathcal{S}_1$ from Eq.(\ref{10y}) instead of $\bar{\mathcal{S}}_1$ but the total cost increases to 7.2944...
\subsubsection{$\boldsymbol{L_{8}(3)}$}For this $L$--function value the following algorithmically different identities are proven from Eq.(\ref{L2}) applying reduction formulas Eq.(\ref{8t}), Eq.(\ref{8v}) and taking $x\in\{\tfrac18,\tfrac38\}$ (for mirror identities $x\in\{\tfrac58,\tfrac78\}$) on the accelerated $\zeta(3,x)$ series $\mathcal{S}_1(x)$ Eq.(\ref{10y})\vspace{-4pt}
\begin{equation}\label{L83a}
	\negthinspace128\,L_8(3)=\ \mathcal{S}_1(\tfrac18)-112\,\zeta(3)-\pi^3\left(4+3\sqrt2\right)\vspace{-6pt}
\end{equation}%
\begin{equation}\label{L83b}
	128\,L_8(3)=-\mathcal{S}_1(\tfrac38)+112\,\zeta(3)-\pi^3\left(4-3\sqrt2\right)\vspace{-2pt}
\end{equation}The total binary splitting cost is $C_s=4.87564...$ (series $\mathcal{S}_1$) + $2.05136...$ ($\zeta(3)$, Eq.(\ref{z2})\thinspace) + $0.36748...$ ($\pi$, Chudnovsky) = 7.2944...
\subsubsection{$\boldsymbol{L_{12}(3)}$}For this $L$--Series, the same approach as $L_8(3)$ but taking $x\in\{\tfrac1{12},\tfrac5{12}\}$ (for mirror identities $x\in\{\tfrac7{12},\tfrac{11}{12}\}$) gives now these algorithmically independent relationships\vspace{-4pt} 
\begin{equation}\label{L123a}
	432\,L_{12}(3)=\ \mathcal{S}_1(\tfrac1{12})-364\,\zeta(3)-2\,\pi^3\left(7+4\sqrt3\right)\vspace{-6pt}
\end{equation}%
\begin{equation}\label{L123b}
	432\,L_{12}(3)=-\mathcal{S}_1(\tfrac5{12})+364\,\zeta(3)+2\,\pi^3\left(7-4\sqrt3\right)\vspace{0pt}
\end{equation}Both formulas have the same total cost $C_s=7.2944...$
\section{Testings\smallskip}\vspace*{4pt}Bench tests were conducted with these new formulas on the y-cruncher platform \cite{YEE0} to determine their real-world performance. This is a specialized software for calculating constants to a huge number of digits limited only by hardware capacity. For hypergeometric type series Eq.(\ref{7b}) y-cruncher's fast computing only accepts rational terms (irrational factors are just applied outside the series) with polynomials $p,q,r$ having coefficients with integer size less than 64--bit\footnote{Bit--size $<$ 64 and minimal cost were the selection conditions to take the most efficient $\zeta(2,x)-\zeta(2,y)$ and $\zeta(3,x)$ accelerated series from WZ--Dougall's $_5H_5$ discovered series output.}.\newline Two custom configuration files for y-cruncher with algorithmically independent formulas were written for each $L$--function value, one is applied for decimal places computation and the other for digits verification. Experimental results are placed in \textit{Tables} 1,2 next page.\vspace{4pt} 

\pagebreak

\subsection{$L_{k}(s)$ Testings}%
\[
\vspace*{4pt}%
\begin{tabular}
	[c]{||l|c|r|r||c|c|r|r||}\hline\hline
	{\small \textit{\small L--Series}} & {\small \textit{\small Equation}} &
	{\small \textit{\small Cost }}$\mathit{C}_{s}$ & {\small \textit{\small Time }[\textit{s}]}
	& {\small \textit{\small L--Series}} & {\small \textit{\small Equation}} &
	{\small \textit{\small Cost }}$\mathit{C}_{s}$ & {\small \textit{\small Time }[\textit{s}%
		]}\\\hline\hline
	\multicolumn{1}{||c|}{${\small L}_{-7}${\small (2)}} & {\small Eq.(\ref{L72})} &
	{\small 8.332} & {\small 338.95} & ${\small L}_{-24}${\small (2)} &
	\multicolumn{1}{|c|}{{\small Eq.(\ref{L24a})}} & {\small 8.292} & {\small 327.43}\\
	\multicolumn{1}{||c|}{} & {\small Eq.(\ref{L72a})} & {\small 14.627} & {\small 555.80}
	&  & \multicolumn{1}{|c|}{{\small Eq.(\ref{L24b})}} & {\small 8.292} & {\small 325.17}%
	\\\hline
	\multicolumn{1}{||c|}{${\small L}_{-8}${\small (2)}} & {\small Eq.(\ref{L82c})} &
	{\small 4.978} & {\small 197.55} & ${\small L}_{5}${\small (3)} &
	\multicolumn{1}{|c|}{{\small Eq.(\ref{L53a})}} & {\small 7.241} & {\small 266.56}\\
	\multicolumn{1}{||c|}{} & {\small Eq.(\ref{L82d})} & {\small 4.978} & {\small 198.92} &
	& \multicolumn{1}{|c|}{{\small Eq.(\ref{L53b})}} & {\small 7.241} & {\small 267.08}%
	\\\hline
	\multicolumn{1}{||c|}{${\small L}_{-15}${\small (2)}} & {\small Eq.(\ref{L15a})} &
	{\small 8.292} & {\small 298.47} & ${\small L}_{8}${\small (3)} &
	\multicolumn{1}{|c|}{{\small Eq.(\ref{L83a})}} & {\small 7.294} & {\small 259.51}\\
	\multicolumn{1}{||c|}{} & {\small Eq.(\ref{L15c})} & {\small 8.292} & {\small 314.49} &
	& \multicolumn{1}{|c|}{{\small Eq.(\ref{L83b})}} & {\small 7.294} & {\small 266.45}%
	\\\hline
	\multicolumn{1}{||c|}{${\small L}_{-20}${\small (2)}} & {\small Eq.(\ref{L20a})} &
	{\small 8.058} & {\small 307.70} & ${\small L}_{12}${\small (3)} &
	{\small Eq.(\ref{L123a})} & {\small 7.294} & {\small 258.76}\\
	\multicolumn{1}{||c|}{} & {\small Eq.(\ref{L20c})} & {\small 8.058} & {\small 316.98} &
	& {\small Eq.(\ref{L123b})} & {\small 7.294} & {\small 264.48}\\\hline\hline
	\multicolumn{8}{||c||}{$%
		\genfrac{}{}{0pt}{0}{\genfrac{}{}{0pt}{1}{{}}{{}}}{{}}%
		$\text{\small Table 1. }\textit{\small Dirichlet L-Series Algorithms Performance}$%
		\genfrac{}{}{0pt}{0}{\genfrac{}{}{0pt}{1}{{}}{{}}}{{}}%
		$}\\
	\multicolumn{8}{||c||}{{\small Time to get 100 million decimal
			digits}}\\\hline\hline
\end{tabular}
\
\]\vspace{-12pt}
\begin{center}%
	\begin{tabular}
		[c]{||lcc||lcc||}\hline\hline
		\multicolumn{6}{||c||}{{\small \text{Table}} {\footnotesize 2.
			}{\small \ \textit{L}}$_{k}${\small (\textit{s})  {\footnotesize Decimal
					Places Verification} }${\small \vspace*{1pt}}$}\\\hline\hline
		{\small \textit{L}}$_{-7}${\small (}{\footnotesize 2}{\small )} & ${\small =}$
		& {\footnotesize 1.1519254705 . . . 8405729486} & {\small \textit{L}}$_{-24}%
		${\small (}{\footnotesize 2}{\small )} & ${\small =\vspace*{1pt}}$ &
		{\footnotesize 1.0578066132 . . . 3165269170}\\
		{\small \textit{L}}$_{-8}${\small (}{\footnotesize 2}{\small )} &
		${\small =\vspace*{1pt}}$ & {\footnotesize 1.0647341710 . . . 9752063771} &
		{\small \textit{L}}$_{5}${\small (}{\footnotesize 3}{\small )} & ${\small =}$
		& {\footnotesize 0.8548247666 . . . 5218792146}\\
		{\small \textit{L}}$_{-15}${\small (}{\footnotesize 2}{\small )} &
		${\small =}$ & {\footnotesize 1.2966185966 . . . 7257465321} & {\small \textit{L}%
		}$_{8}${\small (}{\footnotesize 3}{\small )} & ${\small =\vspace*{1pt}}$ &
		{\footnotesize 0.9583804545 . . . 9019730859}\\
		{\small \textit{L}}$_{-20}${\small (}{\footnotesize 2}{\small )} &
		${\small =\vspace*{1pt}}$ & {\footnotesize 1.1280433247 . . . 0577195474} &
		{\small \textit{L}}$_{12}${\small (}{\footnotesize 3}{\small )} & ${\small =}$
		& {\footnotesize 0.9900400194 . . . 3823927396}\\\hline\hline
		\multicolumn{6}{||c||}{{\footnotesize Last digits from position  99,999,991
				to}{\small \textit{ }}{\footnotesize 100,000,000}${\small \vspace*{1pt}}$%
		}\\\hline\hline
	\end{tabular}
\end{center}\vspace{4pt}
Each $L$--Series was calculated (for first time) up to 100 million decimal digits and elapsed time was measured. A standard 2024 laptop running 64--bit MS Windows\texttrademark\hspace{0pt} 11 OS with no special hardware was used. Timings are the average of five runs. Even for the worst case, $L_{-7}(2)$ Eq.(\ref{L72a}), computing takes under 10 minutes.

\subsection{y--cruncher's custom configuration files}All 16 $L$--Series formulas were written as y-cruncher's configuration (input) files and located embedded as a large comment toward the end of the source TeX file. This TeX file shall be downloaded to extract (copy) the corresponding comment section, paste it into a text editor and split it to save the files. Directions are placed at the 
beginning of such comment. The interested reader can install y-cruncher software \cite{YEE0} and apply these algorithms with these ready to use input files to increase the number of known digits beyond the current limit of $10^8$ decimal places.\vspace*{0pt}

\section{Comments and Conclusions\medskip}
By using the $L$--function definition of Dirichlet $L$--Series depending on Hurwitz $\zeta$'s and applying reduction formulas, several new identities for selected characters that give $L_k(s)$ in terms of $\zeta(2,x)-\zeta(2,y)$ for $s=2$ or $\zeta(3,x)$ for $s=3$ and some additional constants were derived. By applying WZ--Dougall's $_5H_5$ method these (differences of) Hurwitz $\zeta$'s are replaced by accelerated series producing the currently fastest known formulas to calculate such Dirichlet $L$--function values. The new proven formulas were tested to compute these constants up to several million decimal places in a reasonable short time.     

\vspace*{0pt}\medskip

\subsubsection{Acknowledgement}
Jes\'us Guillera has passed away few days ago after a long illness. It is a deeply felt loss. My sincere appreciation to him for his kind support, the original idea of accelerating Hurwitz via WZ to get efficient $L$--Series is his \cite{GUI}. My heartfelt gratitude also for his visionary concepts and several talks that encouraged me to write this article presenting these new proven formulas together with their experimental results.\

\section{Appendix}\vspace{8pt}

\begin{center}
	{\footnotesize \hspace{15pt}ASSOCIATED LMFDB \cite{LMFDB} $L$--FUNCTION DATA LABELS}\vspace{-4pt}
\end{center}
\begin{center}%
	\begin{tabular}
		[c]{lllllll}
		&  &  &  &  &  & \\
		{\small \textit{L}}$_{-3}${\small (}{\footnotesize 2}{\small )} &
		$\rightarrow$ & {\footnotesize 1--3--3.2--{\small \textit{r}}1--0--0} &
		\ \ \  & {\small \textit{L}}$_{-20}${\small (}{\footnotesize 2}{\small )} &
		$\rightarrow$ & {\footnotesize 1--20--20.19--{\small \textit{r}}1--0--0}\\
		{\small \textit{L}}$_{-4}${\small (}{\footnotesize 2}{\small )} &
		$\rightarrow$ & {\footnotesize 1--2{\small \textit{e}}%
			2--4.3--{\small \textit{r}}1--0--0} &  \ \ \  & {\small \textit{L}}$_{-24}%
		${\small (}{\footnotesize 2}{\small )} & $\rightarrow$ &
		{\footnotesize 1--24--24.5--{\small \textit{r}}1--0--0}\\
		{\small \textit{L}}$_{-7}${\small (}{\footnotesize 2}{\small )} &
		$\rightarrow$ & {\footnotesize 1--7--7.6--{\small \textit{r}}1--0--0} &
		\ \ \  & {\small \textit{L}}$_{5}${\small (}{\footnotesize 3}{\small )} &
		$\rightarrow$ & {\footnotesize 1--5--5.4--{\small \textit{r}}0--0--0}\\
		{\small \textit{L}}$_{-8}${\small (}{\footnotesize 2}{\small )} &
		$\rightarrow$ & {\footnotesize 1--2{\small \textit{e}}%
			3--8.3--{\small \textit{r}}1--0--0} &  \ \ \  & {\small \textit{L}}$_{8}%
		${\small (}{\footnotesize 3}{\small )} & $\rightarrow$ &
		{\footnotesize 1--2{\small \textit{e}}3--8.5--{\small \textit{r}}0--0--0}\\
		{\small \textit{L}}$_{-15}${\small (}{\footnotesize 2}{\small )} &
		$\rightarrow$ & {\footnotesize 1--15--15.14--{\small \textit{r}}1--0--0} &
		\ \ \  & {\small \textit{L}}$_{12}${\small (}{\footnotesize 3}{\small )} &
		$\rightarrow$ & {\footnotesize 1--12--12.11--{\small \textit{r}}0--0--0}%
	\end{tabular}
\end{center}\vspace{8pt}
\begin{center}
	{\footnotesize \hspace{45pt}DIRICHLET $L$--SERIES REDUCED IDENTITIES\newline FOR WZ--DOUGALL'S $_5H_5$ ACCELERATED FORMULAS}\vspace{8pt}
\end{center}
\begin{center}
	\
	\begin{tabular}
		[c]{ccl}%
		{\footnotesize 32\thinspace}{\small \textit{L}}$_{-8}${\small (2)} & $=\vspace{1pt}$ &
		$\zeta${\small ({\footnotesize 2,\thinspace3/8}) }$-\ \zeta${\small (}%
		{\footnotesize 2,\thinspace7/8}{\small ) }$+\ ${\footnotesize 2\thinspace}%
		$\pi^{2}\sqrt{\text{{\footnotesize 2}}}$\\
		& $=\vspace{1pt}$ & $\zeta${\small ({\footnotesize 2,\thinspace1/8}) }$-\ \zeta$%
		{\small (}{\footnotesize 2,\thinspace5/8}{\small ) }$-\ ${\footnotesize 2}%
		{\small \thinspace}$\pi^{2}\sqrt{\text{{\footnotesize 2}}}$\\
		{\footnotesize 225\thinspace}{\small \textit{L}}$_{-15}${\small (2)} & $=\vspace{1pt}$ &
		{\footnotesize 4\thinspace}$\zeta${\small ({\footnotesize 2,\thinspace1/15})
		}$-$ {\footnotesize 4\thinspace}$\zeta${\small (}{\footnotesize 2,\thinspace
			11/15}{\small ) }$-$ {\footnotesize 234\thinspace}{\small \textit{L}}$_{-3}%
		${\small (2)}{\footnotesize \thinspace}$-\ ${\small 8\thinspace}$\pi^{2}%
		\sqrt{\text{{\footnotesize 15}}+\text{{\footnotesize 6}}\sqrt
			{\text{{\footnotesize 5}}}}$\\
		& $=\vspace{1pt}$ & {\footnotesize 4\thinspace}$\zeta${\small ({\footnotesize 2,\thinspace
				4/15}) }$-$ {\footnotesize 4\thinspace}$\zeta$%
		{\small ({\footnotesize 2,\thinspace14/15}) }$-$ {\footnotesize 234\thinspace
		}{\small \textit{L}}$_{-3}${\small (2)}{\footnotesize \thinspace}%
		$+\ ${\small 8\thinspace}$\pi^{2}\sqrt{\text{{\footnotesize 15}}%
			+\text{{\footnotesize 6}}\sqrt{\text{{\footnotesize 5}}}}$\\
		& $=\vspace{1pt}$ & {\footnotesize 4\thinspace}$\zeta${\small ({\footnotesize 2,\thinspace
				2/15}) }$-$ {\footnotesize 4\thinspace}$\zeta$%
		{\small ({\footnotesize 2,\thinspace7/15}) }$+$ {\footnotesize 234\thinspace
		}{\small \textit{L}}$_{-3}${\small (2)}{\footnotesize \thinspace}%
		$-\ ${\small 8\thinspace}$\pi^{2}\sqrt{\text{{\footnotesize 15}}%
			-\text{{\footnotesize 6}}\sqrt{\text{{\footnotesize 5}}}}$\\
		& $=\vspace{1pt}$ & {\footnotesize 4\thinspace}$\zeta${\small ({\footnotesize 2,\thinspace
				8/15}) }$-$ {\footnotesize 4\thinspace}$\zeta$%
		{\small ({\footnotesize 2,\thinspace13/15}) }$+$ {\footnotesize 234\thinspace
		}{\small \textit{L}}$_{-3}${\small (2)}{\footnotesize \thinspace}%
		$+\ ${\small 8\thinspace}$\pi^{2}\sqrt{\text{{\footnotesize 15}}%
			-\text{{\footnotesize 6}}\sqrt{\text{{\footnotesize 5}}}}$\\
		{\footnotesize 100\thinspace}{\small \textit{L}}$_{-20}${\small (2)} & $=\vspace{1pt}$ &
		{\footnotesize \thinspace}$\zeta${\small ({\footnotesize 2,\thinspace1/20})
		}$-$ {\footnotesize \thinspace}$\zeta${\small (}{\footnotesize 2,\thinspace
			11/20}{\small ) }$-$ {\footnotesize 96\thinspace{\small \textit{G}}\thinspace
		}$-\ ${\small 2\thinspace}$\pi^{2}\sqrt{\text{{\footnotesize 50}%
			}+\text{{\footnotesize 22}}\sqrt{\text{{\footnotesize 5}}}}$\\
		& $=\vspace{1pt}$ & {\footnotesize \thinspace}$\zeta${\small ({\footnotesize 2,\thinspace
				9/20}) }$-$ {\footnotesize \thinspace}$\zeta${\small (}%
		{\footnotesize 2,\thinspace19/20}{\small ) }$-$ {\footnotesize 96\thinspace
			{\small \textit{G}}\thinspace}$+\ ${\small 2\thinspace}$\pi^{2}\sqrt
		{\text{{\footnotesize 50}}+\text{{\footnotesize 22}}\sqrt
			{\text{{\footnotesize 5}}}}$\\
		& $=\vspace{1pt}$ & {\footnotesize \thinspace}$\zeta${\small ({\footnotesize 2,\thinspace
				3/20}) }$-$ {\footnotesize \thinspace}$\zeta${\small (}%
		{\footnotesize 2,\thinspace13/20}{\small ) }$+$ {\footnotesize 96\thinspace
			{\small \textit{G}}\thinspace}$-\ ${\small 2\thinspace}$\pi^{2}\sqrt
		{\text{{\footnotesize 50}}-\text{{\footnotesize 22}}\sqrt
			{\text{{\footnotesize 5}}}}$\\
		& $=\vspace{1pt}$ & {\footnotesize \thinspace}$\zeta${\small ({\footnotesize 2,\thinspace
				7/20}) }$-$ {\footnotesize \thinspace}$\zeta${\small (}%
		{\footnotesize 2,\thinspace17/20}{\small ) }$+$ {\footnotesize 96\thinspace
			{\small \textit{G}}\thinspace}$+\ ${\small 2\thinspace}$\pi^{2}\sqrt
		{\text{{\footnotesize 50}}-\text{{\footnotesize 22}}\sqrt
			{\text{{\footnotesize 5}}}}$\\
		{\footnotesize 144\thinspace}{\small \textit{L}}$_{-24}${\small (2)} & $=\vspace{1pt}$ &
		{\footnotesize \thinspace}$\zeta${\small ({\footnotesize 2,\thinspace1/24})
		}$-$ {\footnotesize \thinspace}$\zeta${\small (}{\footnotesize 2,\thinspace
			17/24}{\small ) }$-$ {\footnotesize 180\thinspace}{\small \textit{L}}$_{-3}%
		${\small (2)}\thinspace{\footnotesize \thinspace}$-\ ${\small 2\thinspace}%
		$\pi^{2}(\text{{\footnotesize 4}}\sqrt{\text{{\footnotesize 3}}}%
		\,+${\scriptsize \thinspace}{\footnotesize 3}$\sqrt{\text{{\footnotesize 6}}%
		})$\\
		& $=\vspace{1pt}$ & {\footnotesize \thinspace}$\zeta${\small ({\footnotesize 2,\thinspace
				7/24}) }$-$ {\footnotesize \thinspace}$\zeta${\small (}%
		{\footnotesize 2,\thinspace23/24}{\small ) }$-$ {\footnotesize 180\thinspace
		}{\small \textit{L}}$_{-3}${\small (2)}\thinspace{\footnotesize \thinspace
		}$+\ ${\small 2\thinspace}$\pi^{2}(\text{{\footnotesize 4}}\sqrt
		{\text{{\footnotesize 3}}}\,+${\scriptsize \thinspace}{\footnotesize 3}%
		$\sqrt{\text{{\footnotesize 6}}})$\\
		& $=\vspace{1pt}$ & {\footnotesize \thinspace}$\zeta${\small ({\footnotesize 2,\thinspace
				5/24}) }$-$ {\footnotesize \thinspace}$\zeta${\small (}%
		{\footnotesize 2,\thinspace13/24}{\small ) }$+$ {\footnotesize 180\thinspace
		}{\small \textit{L}}$_{-3}${\small (2)}\thinspace{\footnotesize \thinspace
		}$+\ ${\small 2\thinspace}$\pi^{2}(\text{{\footnotesize 4}}\sqrt
		{\text{{\footnotesize 3}}}\,-${\scriptsize \thinspace}{\footnotesize 3}%
		$\sqrt{\text{{\footnotesize 6}}})$\\
		& $=\vspace{1pt}$ & {\footnotesize \thinspace}$\zeta${\small ({\footnotesize 2,\thinspace
				11/24}) }$-$ {\footnotesize \thinspace}$\zeta${\small (}%
		{\footnotesize 2,\thinspace19/24}{\small ) }$+$ {\footnotesize 180\thinspace
		}{\small \textit{L}}$_{-3}${\small (2)}\thinspace{\footnotesize \thinspace
		}$-\ ${\small 2\thinspace}$\pi^{2}(\text{{\footnotesize 4}}\sqrt
		{\text{{\footnotesize 3}}}\,-${\scriptsize \thinspace}{\footnotesize 3}%
		$\sqrt{\text{{\footnotesize 6}}})$\\
		& $=\vspace{1pt}$ & {\footnotesize \thinspace}$\zeta${\small ({\footnotesize 2,\thinspace
				5/24}) }$-$ {\footnotesize \thinspace}$\zeta${\small (}%
		{\footnotesize 2,\thinspace23/24}{\small ) }$-$ {\footnotesize 160\thinspace
			{\small \textit{G}}}\thinspace{\footnotesize \thinspace}$+\ $%
		{\small 2\thinspace}$\pi^{2}(${\footnotesize 5}$\sqrt{\text{{\footnotesize 2}%
		}}\,+${\scriptsize \thinspace}{\footnotesize 4}$\sqrt{\text{{\footnotesize 3}%
		}})$\\
		& $=\vspace{1pt}$ & {\footnotesize \thinspace}$\zeta${\small ({\footnotesize 2,\thinspace
				1/24}) }$-$ {\footnotesize \thinspace}$\zeta${\small (}%
		{\footnotesize 2,\thinspace19/24}{\small ) }$-$ {\footnotesize 160\thinspace
			{\small \textit{G}}}\thinspace{\footnotesize \thinspace}$-\ $%
		{\small 2\thinspace}$\pi^{2}(${\footnotesize 5}$\sqrt{\text{{\footnotesize 2}%
		}}\,+${\scriptsize \thinspace}{\footnotesize 4}$\sqrt{\text{{\footnotesize 3}%
		}})$\\
		& $=\vspace{1pt}$ & {\footnotesize \thinspace}$\zeta${\small ({\footnotesize 2,\thinspace
				11/24}) }$-$ {\footnotesize \thinspace}$\zeta${\small (}%
		{\footnotesize 2,\thinspace17/24}{\small ) }$+$ {\footnotesize 160\thinspace
			{\small \textit{G}}}\thinspace{\footnotesize \thinspace}$+\ $%
		{\small 2\thinspace}$\pi^{2}(${\footnotesize 5}$\sqrt{\text{{\footnotesize 2}%
		}}\,-${\scriptsize \thinspace}{\footnotesize 4}$\sqrt{\text{{\footnotesize 3}%
		}})$\\
		& $=\vspace{1pt}$ & {\footnotesize \thinspace}$\zeta${\small ({\footnotesize 2,\thinspace
				7/24}) }$-$ {\footnotesize \thinspace}$\zeta${\small (}%
		{\footnotesize 2,\thinspace13/24}{\small ) }$+$ {\footnotesize 160\thinspace
			{\small \textit{G}}}\thinspace{\footnotesize \thinspace}$-\ $%
		{\small 2\thinspace}$\pi^{2}(${\footnotesize 5}$\sqrt{\text{{\footnotesize 2}%
		}}\,-${\scriptsize \thinspace}{\footnotesize 4}$\sqrt{\text{{\footnotesize 3}%
		}})$\\
		{\footnotesize 625\thinspace}{\small \textit{L}}$_{5}${\small (3)} & $=\vspace{1pt}$ &
		{\footnotesize 20\thinspace}$\zeta${\small ({\footnotesize 3,\thinspace1/5})
		}$-$ {\footnotesize 620\thinspace}$\zeta${\small ({\footnotesize 3}) }%
		$-\ ${\small 4\thinspace}$\pi^{3}(\sqrt{\text{{\footnotesize 25}%
			}+\text{{\footnotesize 2}}\sqrt{\text{{\footnotesize 5}}}}+\sqrt
		{\text{{\footnotesize 25}}-\text{{\footnotesize 2}}\sqrt
			{\text{{\footnotesize 5}}}\,})$\\
		& $=\vspace{1pt}$ & {\footnotesize 20\thinspace}$\zeta$%
		{\small ({\footnotesize 3,\thinspace4/5}) }$-$ {\footnotesize 620\thinspace
		}$\zeta${\small ({\footnotesize 3}) }$+\ ${\small 4\thinspace}$\pi^{3}%
		(\sqrt{\text{{\footnotesize 25}}+\text{{\footnotesize 2}}\sqrt
			{\text{{\footnotesize 5}}}}+\sqrt{\text{{\footnotesize 25}}%
			-\text{{\footnotesize 2}}\sqrt{\text{{\footnotesize 5}}}}\,)$\\
		& $=\vspace{1pt}$ & $-${\footnotesize 20\thinspace}$\zeta$%
		{\small ({\footnotesize 3,\thinspace2/5}) }$+$ {\footnotesize 620\thinspace
		}$\zeta${\small ({\footnotesize 3}) }$+\ ${\small 4\thinspace}$\pi^{3}%
		(\sqrt{\text{{\footnotesize 25}}+\text{{\footnotesize 2}}\sqrt
			{\text{{\footnotesize 5}}}}-\sqrt{\text{{\footnotesize 25}}%
			-\text{{\footnotesize 2}}\sqrt{\text{{\footnotesize 5}}}}\,)$\\
		& $=\vspace{1pt}$ & $-${\footnotesize 20\thinspace}$\zeta$%
		{\small ({\footnotesize 3,\thinspace3/5}) }$+$ {\footnotesize 620\thinspace
		}$\zeta${\small ({\footnotesize 3}) }$-\ ${\small 4\thinspace}$\pi^{3}%
		(\sqrt{\text{{\footnotesize 25}}+\text{{\footnotesize 2}}\sqrt
			{\text{{\footnotesize 5}}}}-\sqrt{\text{{\footnotesize 25}}%
			-\text{{\footnotesize 2}}\sqrt{\text{{\footnotesize 5}}}\,})$\\
		{\footnotesize 128\thinspace}{\small \textit{L}}$_{8}${\small (3)} & $=\vspace{1pt}$ &
		$\zeta${\small ({\footnotesize 3,\thinspace1/8}) }$-$
		{\footnotesize 112\thinspace}$\zeta${\small ({\footnotesize 3}) }%
		$-\ ${\small \thinspace}$\pi^{3}(${\footnotesize 4}${\small \,}+\,$%
		{\footnotesize 3}$\sqrt{\text{{\footnotesize 2}}})$\\
		& $=\vspace{1pt}$ & $-$ $\zeta${\small ({\footnotesize 3,\thinspace3/8}) }$+$
		{\footnotesize 112\thinspace}$\zeta${\small ({\footnotesize 3}) }%
		$-\ ${\small \thinspace}$\pi^{3}(${\footnotesize 4}${\small \,}-\,$%
		{\footnotesize 3}$\sqrt{\text{{\footnotesize 2}}})$\\
		& $=\vspace{1pt}$ & $-$ $\zeta${\small ({\footnotesize 3,\thinspace5/8}) }$+$
		{\footnotesize 112\thinspace}$\zeta${\small ({\footnotesize 3}) }%
		$+\ ${\small \thinspace}$\pi^{3}(${\footnotesize 4}${\small \,}-\,$%
		{\footnotesize 3}$\sqrt{\text{{\footnotesize 2}}})$\\
		& $=\vspace{1pt}$ & $\zeta${\small ({\footnotesize 3,\thinspace7/8}) }$-$
		{\footnotesize 112\thinspace}$\zeta${\small ({\footnotesize 3}) }%
		$+\ ${\small \thinspace}$\pi^{3}(${\footnotesize 4}${\small \,}+\,$%
		{\footnotesize 3}$\sqrt{\text{{\footnotesize 2}}})$\\
		{\footnotesize 432\thinspace}{\small \textit{L}}$_{12}${\small (3)} & $=\vspace{1pt}$ &
		$\zeta${\small ({\footnotesize 3,\thinspace1/12}) }$-$
		{\footnotesize 364\thinspace}$\zeta${\small ({\footnotesize 3}) }%
		$-\ ${\footnotesize 2}{\small \thinspace}$\pi^{3}(${\footnotesize 7}%
		${\small \,}+\,${\footnotesize 4}$\sqrt{\text{{\footnotesize 3}}})$\\
		& $=\vspace{1pt}$ & $-$ $\zeta${\small ({\footnotesize 3,\thinspace5/12}) }$+$
		{\footnotesize 364\thinspace}$\zeta${\small ({\footnotesize 3}) }%
		$+\ ${\footnotesize 2}{\small \thinspace}$\pi^{3}(${\footnotesize 7}%
		${\small \,}-\,${\footnotesize 4}$\sqrt{\text{{\footnotesize 3}}})$\\
		& $=\vspace{1pt}$ & $-$ $\zeta${\small ({\footnotesize 3,\thinspace7/12}) }$+$
		{\footnotesize 364\thinspace}$\zeta${\small ({\footnotesize 3}) }%
		$-\ ${\footnotesize 2}{\small \thinspace}$\pi^{3}(${\footnotesize 7}%
		${\small \,}-\,${\footnotesize 4}$\sqrt{\text{{\footnotesize 3}}})$\\
		& $=\vspace{1pt}$ & $\zeta${\small ({\footnotesize 3,\thinspace11/12}) }$-$
		{\footnotesize 364\thinspace}$\zeta${\small ({\footnotesize 3}) }%
		$+\ ${\footnotesize 2}{\small \thinspace}$\pi^{3}(${\footnotesize 7}%
		${\small \,}+\,${\footnotesize 4}$\sqrt{\text{{\footnotesize 3}}})$%
	\end{tabular}
	
\end{center}\vspace{8pt}
\begin{center}
	{\footnotesize \hspace{45pt}WZ--DOUGALL'S $_5H_5$ MAPLE COMMANDS \newline FOR PROVING ACCELERATED $L$--SERIES FORMULAS}\vspace{0pt}
\end{center}\vspace{0pt}
\begin{center}
	\
	\begin{tabular}
		[c]{lcll}%
		{\footnotesize \thinspace}{\small \textit{L}}$_{-3}${\small (2)} &  &
		{\small \texttt{SER(1+3\hspace*{1pt}n,1/3+n,1/3+n,2/3+n,2/3+n)}} & {\footnotesize \cite{GUILLERA}}\\
		{\footnotesize \thinspace}{\small \textit{L}}$_{-4}${\small (2)} &
		{\small \textit{G}} & {\small \texttt{SER(1+3\hspace*{1pt}n,1/4+n,1/4+n,3/4+n,3/4+n)}}
		& {\footnotesize \cite{HESS} Eq.(13)}\\
		&  & {\small \texttt{SER(1+4\hspace*{1pt}n,1/4+n,1/4+n,3/4+n,3/4+n)}} & {\footnotesize \cite{BRAD} Eq.(77)}\\
		&  & {\small \texttt{SER(1+8\hspace*{1pt}n,1/2+4\hspace*{1pt}n,1/2+3\hspace*{1pt}n,1/2+2\hspace*{1pt}n,oo)}} & {\footnotesize Eq.(\ref{g3})}\\
			{\small \textit{L}}$_{-7}${\small (2)} &  & {\small \texttt{SER(1+5\hspace*{1pt}%
					n,}$\mathtt{\ell}$\texttt{/7+2\hspace*{1pt}n,}$\mathtt{\ell}$%
				\texttt{/7+2\hspace*{1pt}n,1-}$\mathtt{\ell}$/7\texttt{+2\hspace*{1pt}%
					n,1-}$\mathtt{\ell}$/7\texttt{+2\hspace*{1pt}n)}} & {\footnotesize Eq.(\ref{L72a})}\\
			{\small \textit{L}}$_{-8}${\small (2)} &  & {\small \texttt{SER(5/4+4\hspace
					*{1pt}n,3/8+n,3/8+n,7/8+n,7/8+n)}} & {\footnotesize Eq.(\ref{L82c})}\\
			&  & {\small \texttt{SER(3/4+4\hspace*{1pt}n,1/8+n,1/8+n,5/8+n,5/8+n)}} & {\footnotesize Eq.(\ref{L82d})}\\
			{\small \textit{L}}$_{-15}${\small (2)} &  & {\small \texttt{SER(4/5+3\hspace
					*{1pt}n,1/15+n,1/15+n,11/15+n,11/15+n)}} & {\footnotesize Eq.(\ref{L15a})}\\
			&  & {\small \texttt{SER(3/5+3\hspace*{1pt}n,2/15+n,2/15+n,7/15+n,7/15+n)}} &
			{\footnotesize Eq.(\ref{L15c})}\\
			&  & {\small \texttt{SER(6/5+3\hspace*{1pt}n,4/15+n,4/15+n,14/15+n,14/15+n)}} &
			{\footnotesize p.16}\\
			&  & {\small \texttt{SER(7/5+3\hspace*{1pt}n,8/15+n,8/15+n,13/15+n,13/15+n)}} &
			{\footnotesize p.16}\\
			{\small \textit{L}}$_{-20}${\small (2)} &  & {\small \texttt{SER(3/5+3\hspace
					*{1pt}n,1/20+n,1/20+n,11/20+n,11/20+n)}} & {\footnotesize Eq.(\ref{L20a})}\\
			&  & {\small \texttt{SER(4/5+3\hspace*{1pt}n,3/20+n,3/20+n,13/20+n,13/20+n)}} &
			{\footnotesize Eq.(\ref{L20c})}\\
			&  & {\small \texttt{SER(7/5+3\hspace*{1pt}n,9/20+n,9/20+n,19/20+n,19/20+n)}} &
			{\footnotesize p.16}\\
			&  & {\small \texttt{SER(6/5+3\hspace*{1pt}n,7/20+n,7/20+n,17/20+n,17/20+n)}} &
			{\footnotesize p.16}\\
			{\small \textit{L}}$_{-24}${\small (2)} &  & {\small \texttt{SER(3/4+3\hspace
					*{1pt}n,1/24+n,1/24+n,17/24+n,17/24+n)}} & {\footnotesize Eq.(\ref{L24a})}\\
			&  & {\small \texttt{SER(3/4+3\hspace*{1pt}n,5/24+n,5/24+n,13/24+n,13/24+n)}} &
			{\footnotesize Eq.(\ref{L24b})}\\
			&  & {\small \texttt{SER(5/4+3\hspace*{1pt}n,7/24+n,7/24+n,23/24+n,23/24+n)}} &
			{\footnotesize p.16}\\
			&  & {\small \texttt{SER(5/4+3\hspace*{1pt}n,11/24+n,11/24+n,19/24+n,19/24+n)}} &
			{\footnotesize p.16}\\
			&  & {\small \texttt{SER(7/6+3\hspace*{1pt}n,5/24+n,5/24+n,23/24+n,23/24+n)}} &
			{\footnotesize Eq.(\ref{L24c})}\\
			&  & {\small \texttt{SER(7/6+3\hspace*{1pt}n,11/24+n,11/24+n,17/24+n,17/24+n)}} &
			{\footnotesize Eq.(\ref{L24d})}\\
			&  & {\small \texttt{SER(5/6+3\hspace*{1pt}n,1/24+n,1/24+n,19/24+n,19/24+n)}} &
			{\footnotesize p.16}\\
			&  & {\small \texttt{SER(5/6+3\hspace*{1pt}n,7/24+n,7/24+n,13/24+n,13/24+n)}} &
			{\footnotesize p.16}\\
			{\small \textit{L}}$_{1}${\small (3)} & ${\footnotesize \zeta}${\small (3)} &
			{\small \texttt{SER(2+11\hspace*{1pt}n,1+4\hspace*{1pt}n,1+3\hspace*{1pt}%
					n,1+2\hspace*{1pt}n,1+n)}} & {\footnotesize Eq.(\ref{z1})}\\
			&  & {\small \texttt{SER(2+10\hspace*{1pt}n,1+4\hspace*{1pt}n,1+3\hspace
					*{1pt}n,1+2\hspace*{1pt}n,1+n)}} & {\footnotesize Eq.(\ref{z4})}\\
			{\small \textit{L}}$_{5}${\small (3)} &  & {\small \texttt{SER(2/5+6\hspace
					*{1pt}n,1/5+3\hspace*{1pt}n,1/5+2\hspace*{1pt}n,1/5+2\hspace*{1pt}%
					n,1/5+2\hspace*{1pt}n)}} & {\footnotesize Eq.(\ref{L53a})}\\
			&  & {\small \texttt{SER(4/5+6\hspace*{1pt}n,2/5+3\hspace*{1pt}n,2/5+2\hspace
					*{1pt}n,2/5+2\hspace*{1pt}n,2/5+2\hspace*{1pt}n)}} & {\footnotesize Eq.(\ref{L53b})}\\
			&  & {\small \texttt{SER(6/5+6\hspace*{1pt}n,3/5+3\hspace*{1pt}n,3/5+2\hspace
					*{1pt}n,3/5+2\hspace*{1pt}n,3/5+2\hspace*{1pt}n)}} & {\footnotesize p.16}\\
			&  & {\small \texttt{SER(8/5+6\hspace*{1pt}n,4/5+3\hspace*{1pt}n,4/5+2\hspace
					*{1pt}n,4/5+2\hspace*{1pt}n,4/5+2\hspace*{1pt}n)}} & {\footnotesize p.16}\\
			{\small \textit{L}}$_{8}${\small (3)} &  & {\small \texttt{SER(1/4+4\hspace
					*{1pt}n,1/8+2\hspace*{1pt}n,1/8+n,1/8+n,1/8+n)}} & {\footnotesize Eq.(\ref{L83a})}\\
			&  & {\small \texttt{SER(3/4+4\hspace*{1pt}n,3/8+2\hspace*{1pt}n,3/8+n,3/8+n,3/8+n)}} &
			{\footnotesize Eq.(\ref{L83b})}\\
			&  & {\small \texttt{SER(5/4+4\hspace*{1pt}n,5/8+2\hspace*{1pt}n,5/8+n,5/8+n,5/8+n)}} &
			{\footnotesize p.16}\\
			&  & {\small \texttt{SER(7/4+4\hspace*{1pt}n,7/8+2\hspace*{1pt}n,7/8+n,7/8+n,7/8+n)}} &
			{\footnotesize p.16}\\
			{\small \textit{L}}$_{12}${\small (3)} &  & {\small \texttt{SER(1/6+4\hspace
					*{1pt}n,1/12+2\hspace*{1pt}n,1/12+n,1/12+n,1/12+n)}} & {\footnotesize Eq.(\ref{L123a})}\\
			&  & {\small \texttt{SER(5/6+4\hspace*{1pt}n,5/12+2\hspace*{1pt}%
					n,5/12+n,5/12+n,5/12+n)}} & {\footnotesize Eq.(\ref{L123b})}\\
			&  & {\small \texttt{SER(7/6+4\hspace*{1pt}n,7/12+2\hspace*{1pt}%
					n,7/12+n,7/12+n,7/12+n)}} & {\footnotesize p.16}\\
			&  & {\small \texttt{SER(11/6+4\hspace*{1pt}n,11/12+2\hspace*{1pt}%
					n,11/12+n,11/12+n,11/12+n)}} & {\footnotesize p.16}%
		\end{tabular}
		
	\end{center}
\vfill\eject

\end{document}